\documentclass[11pt]{article}
\pdfoutput=1
\usepackage[latin1]{inputenc}
\usepackage[T1]{fontenc}
\usepackage{amsthm,amsmath,amssymb}
\usepackage[numbers,sort&compress]{natbib}
\usepackage{graphicx}
\usepackage{a4wide}

\theoremstyle{plain}
\newtheorem{theorem}{Theorem}
\newtheorem{corollary}{Corollary}
\newtheorem{lemma}{Lemma}
\newtheorem{proposition}{Proposition}

\theoremstyle{definition}
\newtheorem{definition}{Definition}
\newtheorem{example}{Example}
\newtheorem{remark}{Remark}

\newcommand{\Cov}{\mathrm{cov}}
\newcommand{\Cor}{\mathrm{cor}}
\newcommand{\eiphi}{\mathrm{e}^{\mathrm{i}\varphi}}
\newcommand{\dd}{\mathrm{d}}

\newcommand{\EE}{\mathbb{E}}
\newcommand{\NN}{\mathbb{N}}
\newcommand{\PP}{\mathbb{P}}
\newcommand{\RR}{\mathbb{R}}
\newcommand{\RRol}{\overline{\mathbb{R}}}
\newcommand{\usc}{u.\,s.\,c.}
\newcommand{\Xieps}{\Xi_{\oplus\varepsilon}}
\newcommand{\ZZ}{\mathbb{Z}}
\newcommand{\tildeZ}{Z}

\title{Random Marked Sets}

\author{
	Felix Ballani$^1$, Zakhar Kabluchko$^2$, Martin Schlather$^3$\\[5mm]\\
	${}^1$Institute for Stochastics, TU Bergakademie Freiberg,\\ 
	D-09596 Freiberg, Germany,\\
	ballani@math.tu-freiberg.de\\[5mm] 
	${}^2$Institute of Stochastics, University of Ulm,\\
	Helmholtzstr. 18, D-89069 Ulm, Germany,\\
	zakhar.kabluchko@uni-ulm.de\\[5mm]
	${}^3$Institute for Mathematical Stochastics, Georgia Augusta University,\\
	Goldschmidtstr. 7, D-37077 Göttingen, Germany,\\
	schlather@math.uni-goettingen.de
}

\begin{document}
\maketitle

{\bf Abstract.}
We aim to link random fields and marked point processes
and therefore 
introduce a new class of stochastic processes which are defined on a random
set in $\mathbb{R}^d$. Unlike for random fields, the mark
covariance function of a marked random set 
is in general not positive definite. This implies
that in many situations the use of simple geostatistical methods appears to be
questionable. Surprisingly, for a special class of processes
based on Gaussian random fields, we do have positive definiteness for the
corresponding mark covariance function and mark correlation function.

\bigskip
{\bf Classification.} Primary: 60G60, 60G55; secondary: 60G15, 60D05



\section{Introduction}
Quantities measured in space are mostly modelled as so-called
regionalized variables under the implicite assumption
that these quantities can,
in principle, be measured everywhere and that the choice of
sampling points does not depend on the values of these quantities. Based on
this assumption, several geostatistical methods like variogram analysis or
kriging can be applied \cite{ref:Cressie1993}. However, there are two
types of situations where 
this assumption does not hold \cite{ref:SchlatherRibeiroDiggle2004} and hence,
uncritical 
use of geostatistical methods might cause incorrect or meaningless results.

The first type of problems is caused by the investigators themselves by some
kind of preferential sampling \cite{ref:DiggleMenezesSu2008}. For instance,
this happens when data are 
sampled only at places where high values of the variable of interest are
expected.
The second type of problems is intrinsic to the investigated object
itself. An obvious situation is the investigation of individuals, e.\,g., trees
in a 
forest, where interactions among individuals are present. In this particular
situation the theory of marked point processes provides a formal framework for
data analysis
\cite{ref:DiggleRibeiroChristensen2003,ref:IllianPenttinenStoyanStoyan2008}. 

In this paper, we would like to draw the reader's attention to some further, deceptive
situations which belongs to the second type and 
where implicit conditioning has been mostly ignored in literature
\cite{ref:KangasMaltamo,ref:LiskiWestman1997,ref:Wallermanetal2002}.
For instance, the investigation of pesticides in soil
is restricted to cropland and the height of forest litter is restricted to
silvicultural areas. In both cases, an unintended preselection cannot be excluded
since environmental conditions directly influence the kind of 
land use. A further, simple example has motivated this work and 
appears when 
the altitude is predicted by geostatistical methods based
on measurements that are taken above sea level only. 

Such kind of conditioning might be considered as minor, but can
cause major effects, nonetheless. We advise caution  because
of the following facts: 

1. Any characteristic, such as the covariance function or the variogram, has
to be understood as a conditional quantity given measurements can be taken at
certain locations. 

2. In general, neither the covariance function is positive definite nor the
variogram is conditionally negative definite. 

Since Gaussian random fields are rather popular, a
bigger part of this paper deals with the following hypothetical model: 
 the sea level is at 0 and 
 the altitude is given by some (smooth) stationary
Gaussian random field $Z$ with mean $-t$ and variance $1$. 
Then we face the following oddities when inference is based on  measurements
above sea level only:

1.  The theoretical variogram is not conditionally negative definite, in general.

2. A naive definition of the covariance function $C(x,y)$
by \[C(x,y)=\EE[Z(x)Z(y)\mid Z(x)\geq0,Z(y)\geq0]-\overline{m}^2\] 
leads in general to a function that is not positive definite,
for any $\overline{m}\in\RR$.

3. A more suitable definition of the covariance function for the 
altitude
above sea level as the conditional covariance given that $Z(x)\geq0$ and
$Z(y)\geq0$ leads to a function that is 
never 
differentiable unless the field is spatially constant.

4. If $t=0$, the conditional 
covariance function 
is
positive definite, although, no random field exists that is independent
of the sampling locations and that can model the 
altitude above sea level.

\bigskip

Before discussing the above set-up in detail, we will introduce a theoretical framework so that both a meaningful definition of
 second-order characteristics is possible and usual random fields as well as
 marked point processes are included as particular cases. For this reason we
 extend the notion of a random upper semi-continuous (\usc) function (taking
 values in $\RRol=[-\infty,\infty]$) on $\RR^d$ such that the domain
 is a random subset of $\RR^d$. To this end we make use of Matheron's
 \cite{ref:Matheron1969} idea and consider the hypograph \[A_f=\{(x,t)\in
 X\times\RRol:\,t\leq f(x)\},\quad X\subseteq\RR^d\] 
of a function $f:X\rightarrow\RRol$. In fact, $A_f$ is closed if and only if
$f$ is \usc\ on closed $X$, and the mapping $f\mapsto A_f$ is a bijection. 

\bigskip

The paper is organized as follows. In Section \ref{sec:RandomMarkedSets}, we formally introduce the notion of a random marked closed set and discuss some examples. In Section \ref{sec:Characteristics}, we generalise the definition of several characteristics for random fields to random marked sets. We show that, in general, they do not share the same definiteness properties as their random field analogues. In Section \ref{sec:Gaussian} we study Gaussian random fields $Z(x)$ given that $Z(x)$ exceeds a certain threshold $t\in\mathbb{R}$. In Section \ref{sec:Differentiability} some results on the differentiablity of the mark covariance function of random marked sets are given. 
In Section 
\ref{sec:Proofs}, we collect the proofs of the statements of the preceding sections.

%
\section{Random marked closed sets}\label{sec:RandomMarkedSets}
Denote by $\RRol=\RR\cup\{-\infty,+\infty\}$ the extended real line. Let
\begin{align*}
\Phi_{usc}=\{(X,f):\,X\subseteq\RR^d\text{ is closed},\,f:X\rightarrow\RRol\text{ is \usc}\}.
\end{align*}
$\Phi_{usc}$ is isomorphic to the system $\mathcal{U}_{cl}$ of all closed sets $A\subseteq\RR^d\times\RRol$ which satisfy
\begin{align}\label{eqn:umbra}
&\forall x\in\RR^d\,\,\forall t\in\RRol: (x,t)\in A\Rightarrow\{x\}\times[-\infty,t]\subseteq A
\end{align}
by the bijection
\begin{align*}\label{eqn:tau}
\tau:\Phi_{usc}&\rightarrow\mathcal{U}_{cl}\\ \nonumber
(X,f)&\mapsto\{(x,t)\in X\times\RRol:\,t\leq f(x)\} ,\quad (X,f)\in\Phi_{usc}.
\end{align*}
The subsequent proposition follows immediately from the fact that the space $\mathcal{F}(\RR^d\times\RRol)$ of closed subsets of $\RR^d\times\RRol$ is compact \cite{ref:Matheron1969,ref:Molchanov2005} and that $\mathcal{U}_{cl}$ is closed in $\mathcal{F}(\RR^d\times\RRol)$.
%
\begin{proposition}
$\Phi_{usc}$ is compact in the topology induced by $\mathcal{U}_{cl}$.
\end{proposition}
%
\begin{definition}\label{def:RMS2}
Let $(\Omega,\mathcal{A},\PP)$ be a complete probability space and let $(\Xi,Z):\Omega\rightarrow\Phi_{usc}$ be a mapping with
\begin{align*}
\{\omega\in\Omega:\,\tau((\Xi,Z)(\omega))\cap B\neq\varnothing\}\in\mathcal{A}
\end{align*}
for every compact set $B$ in $\RR^d\times\RRol$. Then $(\Xi,Z)$ is called a \emph{random marked closed set}.
\end{definition}
The distribution law of a random closed set is characterized by the probabilities of hitting compact sets \cite{ref:Matheron1975,ref:Molchanov2005} whereas, by \cite[Prop. 2.3.1]{ref:Matheron1975}, it suffices to consider a suitable base. When choosing the same base of all finite unions of half cylinders $B_i\times[t_i,\infty]$ as in \cite[Thm. XII-6]{ref:Serra1982} for random \usc\ functions on $\RR^d$, we obtain the following characterization of random marked closed sets.
%
\begin{theorem}\label{thm:characterization}
The distribution of a random marked closed set $(\Xi,Z)$ (as a probability measure on $\Phi_{usc}$) is completely determined by the joint probabilities
\begin{align*}
\PP\big(\sup_{x\in B_i\cap\Xi}Z(x)<t_i,\,B_i\cap\Xi\neq\varnothing,\,i\in I;\,B_j\cap\Xi=\varnothing,\,j\in\{1,\ldots,n\}\setminus I\big),
\end{align*}
where $B_1,\ldots,B_n$ are compact subsets of $\RR^d$, $t_1,\ldots,t_n\in\RRol$, and $I$ is a subset of $\{1,\ldots,n\}$, $n\in\NN$.
\end{theorem}
%
\begin{definition}\label{def:stationarity}
A random marked closed set $(\Xi,Z)$ is called \emph{stationary} if \[\PP(\tau(\Xi,Z)+(x,0)\in\cdot\,)=\PP(\tau(\Xi,Z)\in\cdot\,)\] for all $x\in\RR^d$, and it is called \emph{isotropic} if \[\PP(\theta\tau(\Xi,Z)\in\cdot\,)=\PP(\tau(\Xi,Z)\in\cdot\,)\] for all rotations $\theta\in SO_{d+1}$ with $\theta(\RR^d\times\{0\})=\RR^d\times\{0\}$.
\end{definition}
%
\begin{example}\label{ex:randomfieldmodel}
A particular model of a random marked closed set that describes an unbiased
sampling of a random field \cite{ref:Takahata1994} is given when
$Z$ is a random \usc\ function on $\RR^d$ that is 
independent of the random closed set
$\Xi$. We 
call $(\Xi,Z)$ a \emph{random-field model}. 

If the data are consistent with a random-field model, any analysis is simplified
considerably since the domain and the marks can be investigated separately
(see also Remark \ref{rem:randomfieldmodel}) by using standard techniques for
random sets \cite{ref:StoyanKendallMecke1995} and for geostatistical data
\cite{ref:Cressie1993,ref:DiggleRibeiroChristensen2003}. For the particular
case of marked point processes, 
 several tests
for the random-field model hypothesis have been developed
\cite{ref:SchlatherRibeiroDiggle2004,ref:GuanShermanCalvin2004}. 
\end{example}
%
\begin{example}\label{ex:distance}
Let $\Xi$ be a random closed set and $Z(x)=d(\partial\Xi,x)$ the Euclidean
distance of $x\in\RR^d$ to the boundary of $\Xi$. Then $Z$ is even continuous
on $\Xi$. Since local maxima of $Z$ are only attained at locations in the
interior of $\Xi$, the random marked set $(\Xi,Z)$ is a random-field model if
and only if $\Xi=\partial\Xi$ almost surely, in which case $Z$ is trivial. 
\end{example}
%
\begin{example}\label{ex:river}
Cressie et al. \cite{ref:CressieFreyHarchSmith2006} consider 
the spatial prediction on a river network. Here, $\Xi$ is the flow of the river
(as a one-dimensional line or a two-dimensional stripe) and $Z$ models the
dissolved oxygen. 
\end{example}
%
\begin{example}\label{ex:hypersurfaces}
Let $\Xi$ be a random closed set represented as a locally finite union of
closed $C^2$-smooth hypersurfaces in $\RR^d$ such that any two hypersurfaces
intersect at most in a set of measure zero with respect to the
$(d-1)$-dimensional Hausdorff measure. For any $x\in\Xi$, the mark $Z(x)$ 
might be
the maximum of the mean curvatures of the hypersurfaces at $x$. 
The mean curvature has its importance, for example, in the analysis of foams
\cite{ref:Kraynik88}.
\end{example}
%
\begin{remark}
In \cite{ref:Molchanov1983} Molchanov studies \emph{labelled random closed sets} in the sense that a random closed set is split into several closed subsets, see also \cite[p. 141]{ref:Molchanov2005}. Here the marks are at the \emph{nominal} scale. Since random marked sets link \emph{real-valued} marked point processes
and \emph{real-valued} random fields, the concept of the present paper may be seen as an implicit generalisation of labelled random closed sets.
\end{remark}
%
\section{Characteristics for random marked closed sets}\label{sec:Characteristics}
For the description of random fields a set of second-order characteristics like the variogram, the covariance function and the correlation function are used \cite{ref:Cressie1993}. In analogy to these summary functions, several second-order characteristics for marked point processes have been introduced as conditional quantities given the existence of points of the respective unmarked point process \cite{ref:Schlather2001,ref:StoyanKendallMecke1995}. Since point processes can be described as random (counting) measures, these quantities have been derived as Radon-Nikodym derivatives of certain second-order moment measures \cite[Section 2.7]{ref:BenesRataj2004}. Nevertheless, random measures are not always appropriate for the definition of second-order characteristics as the following example illustrates.
%
\begin{example}
Let the stationary random closed set $\Xi$ in $\RR^1$ be given by
\[\Xi=\xi+\bigcup_{z\in\ZZ}[2z-p,2z+p]\cup\{2z+1\},\]
where $p\in(0,\tfrac{1}{3})$ and $\xi$ is uniformly distributed on $[0,1]$. Obviously, interpoint distances $r\in(0,2p]$ are only possible if both points belong to the same segment $\xi+[2z-p,2z+p]$, and interpoint distances $r\in(1-p,1+p]$ are only possible if one point belongs to a segment $\xi+[2z-p,2z+p]$ and the other is from one of the singletons, $\{\xi+2z-1\}$ or $\{\xi+2z+1\}$. Since $\PP(0,r\in\Xi)=0$ for all $r\in(1-p,1+p]$, the approach of defining second-order characteristics using a random measure, which is here based on the Lebesgue measure on $\RR^1$, cannot account for segment-singleton point pairs, and hence, these characteristics are undefined for $r\in(1-p,1+p]$. Nonetheless, it does make sense also to consider the correlation of two marks given that the corresponding points are a distance $r$, $r\in(1-p,1+p]$, apart.
\end{example}
In what follows, $B_{\varepsilon}(x)$ denotes the Euclidean ball in $\RR^d$ with centre $x\in\RR^d$ and radius $\varepsilon\geq0$, $\oplus$ denotes Minkowski addition, and we write shortly $\Xieps$ for $\Xi\oplus B_{\varepsilon}(o)$, and $o$ is the origin in $\RR^d$. Furthermore, $\mathbf{1}_A$ denotes the indicator of $A$.

Let $(\Xi,Z)$ be a random marked closed set in $\RR^d$ with marks
in $\RR$.  For ease, we assume stationarity, but the approach can 
be extended to a non-stationary set-up.
For any $\varepsilon\geq0$ define the (stationary) random field $\tildeZ_{\varepsilon}$ by
\begin{align*}
\tildeZ_{\varepsilon}(x)&=\begin{cases}\max\limits_{y\in\Xi\cap B_{\varepsilon}(x)}Z(y),&x\in\Xieps,\\ 0,&\text{otherwise}.\end{cases}
\end{align*}
Let $f:\RR^2\rightarrow\RR$ be a right-continuous function. For all $h\in\RR^d$ define
\begin{align}\label{def:kappa_f}
\kappa_f(h)&=\lim_{\varepsilon\rightarrow0+}\EE\left[f\left(\tildeZ_{\varepsilon}(o),\tildeZ_{\varepsilon}(h)\right)\mid o,h\in\Xieps\right]
\end{align}
whenever $\kappa_{|f|}(h)<\infty$ and $\PP(o,h\in\Xieps)>0$ for all $\varepsilon>0$, otherwise $\kappa_f(h)$ is undefined.

In particular, for the following choices of $f$,
\begin{align}\label{eqn:choices}
e(m_1,m_2)&=m_1,\quad c(m_1,m_2)=m_1m_2,\quad v(m_1,m_2)=m_1^2,
\end{align}
define
\begin{align}
E(h)&=\kappa_e(h)\label{eqn:char:E}\\
\gamma(h)&=\frac{1}{2}(\kappa_v(h)+\kappa_v(-h))-\kappa_c(h)\label{eqn:char:gamma}\\
\Cov(h)&=\kappa_c(h)-\kappa_e(h)\kappa_e(-h)\label{eqn:char:cov}
\\
\Cor(h)&=\frac{\kappa_c(h)-\kappa_e(h)\kappa_e(-h)}{(\kappa_v(h)-\kappa_e(h)^2)^{1/2}(\kappa_v(-h)-\kappa_e(-h)^2)^{1/2}}\label{eqn:char:cor}\\
k_{mm}(h)&=(\overline{m})^{-2}\,\kappa_c(h),\quad(\overline{m}\neq0)\label{eqn:char:kmm},
\end{align}
where $\overline{m}=\mathbb{E}[Z(o)\mid o\in\Xi]$ is the mean mark.

We call $\gamma$ the \emph{mark variogram}, $\mathrm{cov}$ the \emph{mark covariance function}, $\mathrm{cor}$ the \emph{mark correlation function} and $k_{mm}$ \emph{Stoyan's $k_{mm}$-function} of $(\Xi,Z)$ \cite{ref:Schlather2001}. Note that, if $\Xi\equiv\mathbb{R}^d$, these definitions are compatible with the classical definitions for random fields (see Remark \ref{rem:randomfieldmodel}).

Whenever $(\Xi,Z)$ is assumed to be both stationary and isotropic the characteristics given by (\ref{eqn:char:E})--(\ref{eqn:char:kmm}) are rotation invariant. By slight abuse of notation we will write $E(r)$, $r\in[0,\infty)$, instead of $E(h)$, $h\in\RR^d$. The same applies for the functions defined in Eq. (\ref{eqn:char:gamma})--(\ref{eqn:char:kmm}).
%
\begin{remark}
Let $\Psi_{\varepsilon}=\nu_d(\,\cdot\cap\Xieps)$ be the random volume measure associated with the random closed set $\Xieps$. Here, $\nu_d$ is the $d$-dimensional Lebesgue measure. If $\mu^{(2)}_{\varepsilon}$ denotes the second-order moment measure of $\Psi_{\varepsilon}$ then, for $B_1,B_2\in\mathcal{B}(\RR^d)$, we have
\begin{align*}
&\int_{B_2}\int_{B_1}\EE\left[f\left(\tildeZ_{\varepsilon}(x),\tildeZ_{\varepsilon}(y)\right)\mathbf{1}_{\Xieps}(x)\mathbf{1}_{\Xieps)}(y)\right]\,\dd x\,\dd y\\
&\quad=\EE\left[\int_{B_2}\int_{B_1}f\left(\tildeZ_{\varepsilon}(x),\tildeZ_{\varepsilon}(y)\right)\,\Psi_{\varepsilon}(\dd
  x)\,\Psi_{\varepsilon}(\dd y)\right]\\
&\quad=\int_{B_1\times
  B_2}\int_{\RR^2}f(m_1,m_2)\,Q_{\varepsilon;x,y}(\dd(m_1,m_2))\,\mu^{(2)}_{\varepsilon}(\dd(x,y))\\
&\quad=\int_{B_2}\int_{B_1}\int_{\RR^2}f(m_1,m_2)\,Q_{\varepsilon;x,y}(\dd(m_1,m_2))\,\PP(x,y\in\Xieps)\,\dd x\,\dd y,
\end{align*}
where $Q_{\varepsilon;x,y}$ is the two-point mark distribution of the weighted random measure $(\Psi_{\varepsilon},\tildeZ_{\varepsilon})$ \cite{ref:BenesRataj2004}. Hence, for almost all $(x,y)$ with $\PP(x,y\in\Xieps)>0$, we have
\begin{align*}
\EE\left[f\left(\tildeZ_{\varepsilon}(x),\tildeZ_{\varepsilon}(y)\right)\mid x,y\in\Xieps\right]&=\int_{\RR^2}f(m_1,m_2)\,Q_{\varepsilon;x,y}(\dd(m_1,m_2)).
\end{align*}
\end{remark}
%
\begin{remark}\label{rem:simplification}
In case $\PP(o,h\in\Xi)>0$, $h\in\RR^d$, the above definition takes the simpler form
\begin{align}\label{eqn:integrability}
\kappa_f(h)=\EE[f(Z(o),Z(h))\mid o,h\in\Xi]
\end{align}
for $f\in\{e, c, v\}$ if we impose the integrability conditions 
\[
\EE[|Z(o)|^2\mathbf{1}_{\Xi}(o)]<\infty,\;\;
\kappa_{v}(h)<\infty,
\]
and
\[
\lim_{\varepsilon\to0+}\EE[(Z(p(\Xi,o))_-)^2\mathbf{1}_{\Xieps\setminus\Xi}(o)]<\infty.
\]
Here $e$, $c$ and $v$ are given by (\ref{eqn:choices}), $a_-$ denotes the
negative part of $a\in\RR$, and $p(A,x)$ is the metric projection
\cite{ref:HugLastWeil2004} of $x\in\RR^d\setminus\Xi$ onto the
boundary $\partial\Xi$ of $\Xi$ with the smallest coordinates in the
lexicographical ordering, say. Note that the latter is not crucial,
since, due to stationarity, the probability that the projection of $o$
onto $\partial\Xi$ is not unique is zero, see
\cite{ref:HugLastWeil2004}. 

The equality (\ref{eqn:integrability}) can be seen as follows. Denoting
by $a_+$ the positive part of $a\in\RR$ we always
have 
\begin{align*}
Z_{\varepsilon}(o)_+\mathbf{1}_{\Xieps}(o)\mathbf{1}_{\Xieps}(h)&\leq
Z_{\overline{\varepsilon}_1}(o)_+\mathbf{1}_{\Xi_{\oplus\overline{\varepsilon}_1}}(o)\mathbf{1}_{\Xi_{\oplus\overline{\varepsilon}_1}}(h)\\
&\leq
|Z_{\overline{\varepsilon}_1}(o)|\mathbf{1}_{\Xi_{\oplus\overline{\varepsilon}_1}}(o)\mathbf{1}_{\Xi_{\oplus\overline{\varepsilon}_1}}(h)
\end{align*}
for all $0<\varepsilon\leq\overline{\varepsilon}_1$, where the right-hand side
is integrable for $\overline{\varepsilon}_1$ small enough
as  $\kappa_{|e|}(h)<\infty$. Similarly, 
\[Z_{\varepsilon}(o)_-\mathbf{1}_{\Xi}(o)\mathbf{1}_{\Xieps}(h)\leq Z_{\varepsilon}(o)_-\mathbf{1}_{\Xi}(o)\leq Z(o)_-\mathbf{1}_{\Xi}(o)\leq |Z(o)|\mathbf{1}_{\Xi}(o).\]
Finally, since $Z_{\varepsilon}(o)_-=\min_{x\in\Xi\cap B_{\varepsilon}(o)}Z(x)_-$ for $o\in\Xieps\setminus\Xi$, we have
\begin{align*}
Z_{\varepsilon}(o)_-\mathbf{1}_{\Xieps\setminus\Xi}(o)\mathbf{1}_{\Xieps}(h)&\leq Z_{\varepsilon}(o)_-\mathbf{1}_{\Xieps\setminus\Xi}(o)\leq Z(p(\Xi,o))_-\mathbf{1}_{\Xieps\setminus\Xi}(o)\\
&\leq Z(p(\Xi,o))_-\mathbf{1}_{\Xi_{\oplus\overline{\varepsilon}_2}}(o)
\end{align*}
for all $0<\varepsilon\leq\overline{\varepsilon}_2$, where the right-hand side
is integrable for $\overline{\varepsilon}_2$ small enough. Note that
\[
\EE[|Z(o)|\mathbf{1}_{\Xi}(o)]<\infty,\quad\kappa_{|e|}(h)<\infty,\quad
\lim_{\varepsilon\to0+}\EE[Z(p(\Xi,o))_-\mathbf{1}_{\Xieps\setminus\Xi}(o)]<\infty
\]
by Cauchy-Schwarz. 

Similarly, we obtain that
\[|Z_{\overline{\varepsilon}_1}(o)|^2\mathbf{1}_{\Xi_{\oplus\overline{\varepsilon}_1}}(o)\mathbf{1}_{\Xi_{\oplus\overline{\varepsilon}_1}}(h)+|Z(o)|^2\mathbf{1}_{\Xi}(o)+(Z(p(\Xi,o))_-)^2\mathbf{1}_{\Xi_{\oplus\overline{\varepsilon}_2}}(o)\]
is an integrable upper bound of
$|Z_{\varepsilon}(o)|^2\mathbf{1}_{\Xieps}(o)\mathbf{1}_{\Xieps}(h)|$
and is the half part of the
 integrable upper bound of $|Z_{\varepsilon}(o)Z_{\varepsilon}(h)|\mathbf{1}_{\Xieps}(o)\mathbf{1}_{\Xieps}(h)$ due to $|Z_{\varepsilon}(o)Z_{\varepsilon}(h)|\leq\tfrac{1}{2}(|Z_{\varepsilon}(o)|^2+|Z_{\varepsilon}(h)|^2)$. 

Since $Z$ is \usc\ on $\Xi$, a value $\varepsilon>0$ exists for every $x\in\Xi$ and for every $\delta>0$ such that $Z(y)\leq Z(x)+\delta$ for all $y\in B_{\varepsilon}(x)\cap\Xi$. Hence, we have $Z_{\varepsilon}(x)\rightarrow Z(x)$ from above as $\varepsilon\rightarrow0+$. Further, $x\notin\Xi$ implies $x\notin\Xieps$ for all sufficiently small $\varepsilon$. We then have \[f(Z_{\varepsilon}(o),Z_{\varepsilon}(h))\mathbf{1}_{\Xieps}(o)\mathbf{1}_{\Xieps}(h)\rightarrow f(Z(o),Z(h))\mathbf{1}_{\Xi}(o)\mathbf{1}_{\Xi}(h)\text{ a.\,s.}\]
as $\varepsilon\rightarrow0+$. Hence, by the dominated convergence theorem, we have
\begin{align*}
\kappa_f(h)&=\lim_{\varepsilon\rightarrow0+}\frac{\EE\left[f\left(Z_{\varepsilon}(o),Z_{\varepsilon}(h)\right)\mathbf{1}_{\Xieps}(o)\mathbf{1}_{\Xieps}(h)\right]}{\PP(o,h\in\Xieps)}\\
&=\frac{\EE\left[f\left(Z(o),Z(h)\right)\mathbf{1}_{\Xi}(o)\mathbf{1}_{\Xi}(h)\right]}{\PP(o,h\in\Xi)}.
\end{align*}
\end{remark}
%
\begin{remark}
There exists an alternative concept of random marked sets which is inspired by the notion of random fields and where second-order characteristics in the sense of the preceding remark can be defined.

Let $\RRol_{\varnothing}=\RRol\cup\{\zeta_{\varnothing}\}$ be the extension of
$\RRol$ by some $\zeta_{\varnothing}$.  We denote by
$\mathcal{B}(\RRol_{\varnothing})$ the respective Borel $\sigma$-field which
is generated by all sets $B_1\cup B_2$ for $B_1\in\mathcal{B}(\RR)$ and 
$B_2\subset\{-\infty,\infty,\zeta_{\varnothing}\}$.

A family of random variables $Z(\cdot,x):\Omega\rightarrow\RRol_{\varnothing}$, $x\in\RR^d$, on the probability space $(\Omega,\mathcal{A},\PP)$ is called a \emph{random field with random domain} $\Xi$, if
\begin{align*}
\Xi=\{x\in\RR^d:\,Z(\cdot,x)\neq\zeta_{\varnothing}\}.
\end{align*}
Clearly, when $Z$ takes only values different from $\zeta_{\varnothing}$ or $\{-\infty,\infty,\zeta_{\varnothing}\}$ this notion of a random marked set includes usual $\RRol$- or $\RR$-valued random fields on $\RR^d$.

Note that $\Xi$ is a \emph{random set} in a very general sense
\cite{ref:Matheron1975}, entirely determined by its indicator
$\mathbf{1}_{\Xi}(x)=\mathbf{1}_{\RRol}(Z(x))$. If $Z$ is jointly measurable,
i.\,e., $Z$ is
$(\mathcal{A}\otimes\mathcal{B}(\RR^d),\mathcal{B}(\RRol_{\varnothing}))$-measurable,
then the realizations of $\Xi$ are almost surely Borel measurable. If we have
even almost surely closed (open) realizations of $\Xi$ then $Z$ is called a
\emph{random field with random closed (open) domain}, see also
\cite{ref:Molchanov2005}. 

If $\PP(o\in\Xi)>0$ holds for a stationary random field $Z$ with random domain
$\Xi$ we can define second-order characteristics without any further
assumption on path regularity. Let $\widetilde{Z}$ be the (stationary) random
field given by $\widetilde{Z}(x)=Z(x)$ for $x\in\Xi$, and $\widetilde{Z}(x)=0$
otherwise. Let $f:\RR^2\rightarrow\RR$ be a measurable function. For all
$h\in\RR^d$ define 
\begin{align*}
\kappa_f(h)&=\EE[f(\widetilde{Z}(o),\widetilde{Z}(h))\mid o,h\in\Xi]
\end{align*}
whenever $\PP(o,h\in\Xi)>0$ and $\EE[|f(\widetilde{Z}(o),\widetilde{Z}(h))|\mathbf{1}_{\Xi}(o)\mathbf{1}_{\Xi}(h)]<\infty$. 
\end{remark}
%
\begin{remark}\label{rem:randomfieldmodel}
Let $(\Xi,Z)$ be a stationary real-valued random-field model and
\[Z_{\varepsilon}(x)=\begin{cases}\max\limits_{y\in\Xi\cap B_{\varepsilon}(x)}Z(y),&x\in\Xieps,\\ Z(x),&\text{otherwise}.\end{cases}\]
Since $Z$ is \usc\ on $\Xi$ we have $Z_{\varepsilon}(x)\rightarrow Z(x)$ from above for $x\in\Xi$, and hence, by the definition of $Z_{\varepsilon}$, for all $x\in\RR^d$ as $\varepsilon\rightarrow0+$. Then, using the independence of $Z$ and $\Xi$, we obtain
\begin{align*}
\kappa_f(h)&=\lim_{\varepsilon\rightarrow0+}\EE\left[f\left(Z_{\varepsilon}(o),Z_{\varepsilon}(h)\right)\mid o,h\in\Xieps\right]=\lim_{\varepsilon\rightarrow0+}\EE\left[f\left(Z_{\varepsilon}(o),Z_{\varepsilon}(h)\right)\right]\\
&=\EE\left[f\left(Z(o),Z(h)\right)\right]
\end{align*}
for all $h\in\RR^d$ which satisfy $\PP(o,h\in\Xieps)>0$ for all $\varepsilon>0$ and, depending on the choice of $f$ according to (\ref{eqn:choices}), one of the integrability conditions in Remark \ref{rem:simplification} with $\Xi$ replaced by $\RR^d$.
\end{remark}
%
\begin{remark}\label{rem:mpp}
The definition of $\kappa_f$ according to (\ref{def:kappa_f}) is,
in important situations, consistent
with the classical definition of the second-order characteristics of
stationary marked point processes \cite{ref:Schlather2001}. Let
$\widetilde{\Phi}$ be a stationary simple marked point process on
$\RR^d\times\RR$. Then $\Xi$ is the support of the unmarked point process
$\Phi=\widetilde{\Phi}(\,\cdot\times\RR)$. We assume that the second-order
moment measure $\mu^{(2)}$ of $\Phi$ is locally finite. Denoting by $\mathbf{1}_{\{\cdot\}}$ the indicator of the event $\{\cdot\}$, we have 
\begin{align*}
&\EE\left[f\left(\tildeZ_{\varepsilon}(o),\tildeZ_{\varepsilon}(h)\right)\mathbf{1}_{\Xi\oplus B_{\varepsilon}(o)}(o)\mathbf{1}_{\Xi\oplus B_{\varepsilon}(o)}(h)\right]\\
&\quad=\EE\left[f\left(\tildeZ_{\varepsilon}(o),\tildeZ_{\varepsilon}(h)\right)\mathbf{1}_{\{\Phi(B_{\varepsilon}(o))=1\}}\mathbf{1}_{\{\Phi(B_{\varepsilon}(h))=1\}}\right]\\
&\quad\quad+\EE\left[f\left(\tildeZ_{\varepsilon}(o),\tildeZ_{\varepsilon}(h)\right)\mathbf{1}_{\{\Phi(B_{\varepsilon}(o))>1\}}\mathbf{1}_{\{\Phi(B_{\varepsilon}(h))\geq1\}}\right]\\
&\quad\quad+\EE\left[f\left(\tildeZ_{\varepsilon}(o),\tildeZ_{\varepsilon}(h)\right)\mathbf{1}_{\{\Phi(B_{\varepsilon}(o))=1\}}\mathbf{1}_{\{\Phi(B_{\varepsilon}(h))>1\}}\right]
\end{align*}
 for $\|h\|>0$. For any $0<\varepsilon<\|h\|/2$ the first summand equals
\begin{align*}
\EE\left[\sum_{(x_1,m_1),(x_2,m_2)\in\widetilde{\Phi}}f(m_1,m_2)\mathbf{1}_{B_{\varepsilon}(o)}(x_1)\mathbf{1}_{B_{\varepsilon}(h)}(x_2)\right]=:\mu_f^{(2)}(B_{\varepsilon}(o)\times B_{\varepsilon}(h)).
\end{align*}
We can extend the argumentation in \cite[Prop. 9.3.XV]{ref:DaleyVereJones2008} in order to conclude that
\begin{align*}
\frac{\PP(o,h\in\Xieps)}{\mu^{(2)}(B_{\varepsilon}(o)\times B_{\varepsilon}(h))}&=\frac{\PP(\Phi(B_{\varepsilon}(o))\geq1,\Phi(B_{\varepsilon}(h))\geq1)}{\mu^{(2)}(B_{\varepsilon}(o)\times B_{\varepsilon}(h))}\rightarrow1
\end{align*}
as $\varepsilon\rightarrow0+$. 
If we additionally impose the condition that for some
$\overline{\varepsilon}>0$,
{\small
\begin{align*}
\sup_{\varepsilon\in(0,\overline{\varepsilon})}\frac{\EE\left[\left|f\left(\tildeZ_{\varepsilon}(o),\tildeZ_{\varepsilon}(h)\right)\right|\mathbf{1}_{\{\Phi(B_{\varepsilon}(o))>1\}}\mathbf{1}_{\{\Phi(B_{\varepsilon}(h))\geq1\}}\mathbf{1}_{\{|f(\tildeZ_{\varepsilon}(o),\tildeZ_{\varepsilon}(h))|>M\}}\right]}{\PP(\Phi(B_{\varepsilon}(o))\geq1,\Phi(B_{\varepsilon}(h))\geq1)}
\rightarrow 0
\end{align*}
}
as $M\rightarrow\infty$, we obtain
\begin{align*}
\kappa_f(h)=\lim_{\varepsilon\rightarrow0+}\frac{\mu_f^{(2)}(B_{\varepsilon}(o)\times B_{\varepsilon}(h))}{\mu^{(2)}(B_{\varepsilon}(o)\times B_{\varepsilon}(h))}
\end{align*}
which equals $\mu^{(2)}$-a.\,e. the Radon-Nikodym
derivative \[\frac{\mathrm{d}\mu_f^{(2)}(x,x+h)}{\mathrm{d}\mu^{(2)}(x,x+h)}.\]

For instance, the above condition is satisfied if 
$\EE
\left|f\left(\tildeZ_{\varepsilon}(o),\tildeZ_{\varepsilon}(h)\right)\right|^\alpha$
is uniformly bounded on $(0,\overline{\varepsilon})$ for some $\alpha>1$.
\end{remark}

\bigskip

A function $f:\RR^d\rightarrow\RR$ is called \emph{positive definite} if
\[\sum_{i=1}^n\sum_{j=1}^n a_ia_jf(x_i-x_j)\geq0\] for any $n\in\NN$, $x_1,\ldots,x_n\in\RR^d$, and $a_1,\ldots,a_n\in\RR$, and $f$ is called \emph{conditionally negative definite} if \[\sum_{i=1}^n\sum_{j=1}^n a_ia_jf(x_i-x_j)\leq0\]
for any $n\in\NN$, $x_1,\ldots,x_n\in\RR^d$, and for all $a_1,\ldots,a_n\in\RR$ with $\sum_{i=1}^na_i=0$.

For a random-field model all second-order characteristics coincide with those of a random field with \usc\ paths, see Remark \ref{rem:randomfieldmodel}, and thus, share the same definiteness properties. On the other hand, for marked point processes it has been shown by examples \cite{ref:WaelderStoyan1996} and systematically \cite{ref:Schlather2001} that the mark covariance function, the mark correlation function and the $k_{mm}$-function need not be positive definite, and the mark variogram need not be conditionally negative definite in contrast to random fields. Some of the constructions used in \cite{ref:Schlather2001} are based on the fact that for a marked point process, $\Xi$ is a locally finite subset of $\mathbb{R}^d$ and has therefore Lebesgue measure zero. However, the next example shows that in general we cannot expect that the mark covariance function is positive definite (and the mark correlation function and the $k_{mm}$-function either) even when we have $\mathbb{E}[\nu_d(\Xi\cap[0,1]^d)]=P(o\in\Xi)>0$.
%
\begin{example}[Continuation of Example \ref{ex:distance}]\label{ex:periodic}
Let $p\in(\tfrac{2}{3},1]$, $\xi$ be a random variable uniformly distributed on $[0,1]$, $\Xi=\mathbb{Z}\oplus[\xi,p+\xi]$, and $Z(\xi,\cdot\,)$ a 1-periodic function defined by
\begin{align*}
Z(\xi,x)&=\begin{cases}x-\xi, & x\in\mathbb{Z}\oplus[\xi,\tfrac{p}{2}+\xi),\\ p-(x-\xi), & x\in\mathbb{Z}\oplus[\tfrac{p}{2}+\xi,p+\xi),\\ 0, & x\in\mathbb{Z}\oplus[p+\xi,1+\xi).\end{cases}
\end{align*}
Then $Z$ and $\Xi$ are jointly stationary and each of the characteristics given by (\ref{eqn:char:E})--(\ref{eqn:char:kmm}) is 1-periodic. In particular, on $[0,1/2)$ we have 
\begin{align*}
\mathrm{cov}(r)&=
\begin{cases}
\frac{p^4 - 4p^3r - 12p^2r^2 + 48pr^3 - 36r^4}{48(p-r)^2}, & r\in[0,1-p),\\
\frac{ 32r^3(2p-1) + 24r(1-r) - 48pr(1-p)(1-r)   - 12p^2(2pr - p +1) 
-3p^4  +  8p - 4}{48(2p-1)^2}, & r\in[1-p,\tfrac{p}{2}),\\
-\frac{4p^4 - 8p^3 + 6p^2 - 2p + 12r^4 - 24r^3 + 18r^2 - 6r + 1}{12(2p-1)^2}, & r\in[\tfrac{p}{2},\tfrac{1}{2}],
\end{cases}
\end{align*}
and, by symmetry, $\mathrm{cov}(r)=\mathrm{cov}(1-r)$ for $r\in(\tfrac{1}{2},1)$. Since $\mathrm{cov}$ is 1-periodic the 0th coefficient of the Fourier series of $\mathrm{cov}$ is proportional to
\begin{align*}
\int\limits_0^1\mathrm{cov}(r)\,\mathrm{d}r&=
\tfrac{7}{6}p^3\ln\left(\tfrac{p}{2p-1}\right)+\tfrac{409p^5-790p^4+565p^3-280p^2+120p-24}{120(2p-1)^2},
\end{align*}
which is negative for $\frac{2}{3}\leq p<1$ (and vanishes for $p=1$, which is
the random field case). Since $\mathrm{cov}$ is continuous, Bochner's theorem
\cite{ref:Sasvari1994} implies that $\mathrm{cov}$ cannot be a positive
definite function. 
\end{example}
%
\begin{example}\label{ex:covEps}
This example can be seen both as a considerable generalization of Example \ref{ex:distance} and as an attempt to give a set-up for the situation in Example \ref{ex:river}. Let $\Psi$ be a random closed set in $\RR^d$ with almost sure topologically regular realizations and $\Xi$ be the closure of the complement of $\Psi$. Let $Z(x)=f(d(\Psi,x))+Y(x)$, $x\in\RR^d$, where $f$ is any measurable real-valued function of the Euclidean distance between $x$ and $\Psi$, and $Y(x)$ is any real-valued random field on $\RR^d$ independent of $\Psi$. 

In this model, $\Psi$ might, for instance, represent woodland and $\Xi$ cropland. Furthermore, $Z(x)$ might be any quantity of interest concerning cropland, where $Z(x)$ depends, besides additional effects given by $Y(x)$, somehow on the distance to the woodland since, for instance, light and soil properties are influenced by the trees. Accounting for the fact that the influence of trees decreases with increasing distance or vanishes beyond a certain distance, reasonable choices for $f$ are, for instance, $f(t)=c\left(1-\mathrm{e}^{-\alpha t}\right)$, $t\geq0$, $c>0$, $\alpha>0$, or, $f(t)=c\,\max\{t/R,1\}$, $t\geq0$, $c>0$, $R>0$. 

In order to illustrate the problem with the definiteness properties, let $\Psi$ be a stationary Boolean model of intensity $\lambda$ and with compact typical grain $\Psi_0$ \cite{ref:StoyanKendallMecke1995}, $Y\equiv0$, and $f(t)=t$, $t\geq0$, in which case $Z(x)$ is simply the field of contact distances with respect to $\Psi$. Then, for all $s,t\geq0$ and all $h\in\RR^d$, we have
\begin{align*}
\mathbb{P}(d(\Psi,o)>s,d(\Psi,h)>t\mid o,h\in\Xi)=\frac{\exp\{-\lambda u(\Psi_0,h,s,t)\}}{\exp\{-\lambda u(\Psi_0,h,0,0)\}},
\end{align*}
where $u(\Psi_0,h,s,t)=\mathbb{E}[\nu_d((\Psi_0\oplus B_s(o))\cup(\Psi_0\oplus B_t(-h)))]$. Hence, the mark covariance function $\mathrm{cov}$ of $(\Xi,Z)$ is given by 
\begin{align*}
\mathrm{cov}(h)
&=\mathbb{E}[d(\Psi,o)d(\Psi,h)\mid o,h\in\Xi]-(\mathbb{E}[d(\Psi,o)\mid o,h\in\Xi])^2\\
&=\frac{\int_0^{\infty}\int_0^{\infty}\exp\{-\lambda u(\Psi_0,h,s,t)\}\,\mathrm{d}t\,\mathrm{d}s}{\exp\{-\lambda u(\Psi_0,h,0,0)\}}-\frac{\left(\int_0^{\infty}\exp\{-\lambda u(\Psi_0,h,s,0)\}\,\mathrm{d}s\right)^2}{\exp\{-2\lambda u(\Psi_0,h,0,0)\}}.
\end{align*}
For ease, let $d=1$, $\varepsilon\in[0,\infty)$ and $\Psi_0=[-\varepsilon,\varepsilon]$. Then, for $\varepsilon=0$, $\Xi=\RR$, and hence, $\mathrm{cov}$ is positive definite. However, this is not true for any $\varepsilon>0$. Rather than presenting the somewhat lengthy expression for $\mathrm{cov}(r)$, $r\geq0$, which can be given in closed form, we refer the reader to Figure \ref{fig:covEps}. There, the Fourier transform $\widetilde{\mathrm{cov}}$ of $\mathrm{cov}$ indicates that, due to the occurrence of negative values of $\widetilde{\mathrm{cov}}$, $\mathrm{cov}(r)$ is not positive definite for $\varepsilon=1$ and $\varepsilon=2$.
%
\begin{figure}
\begin{center}
\includegraphics[height=6cm]{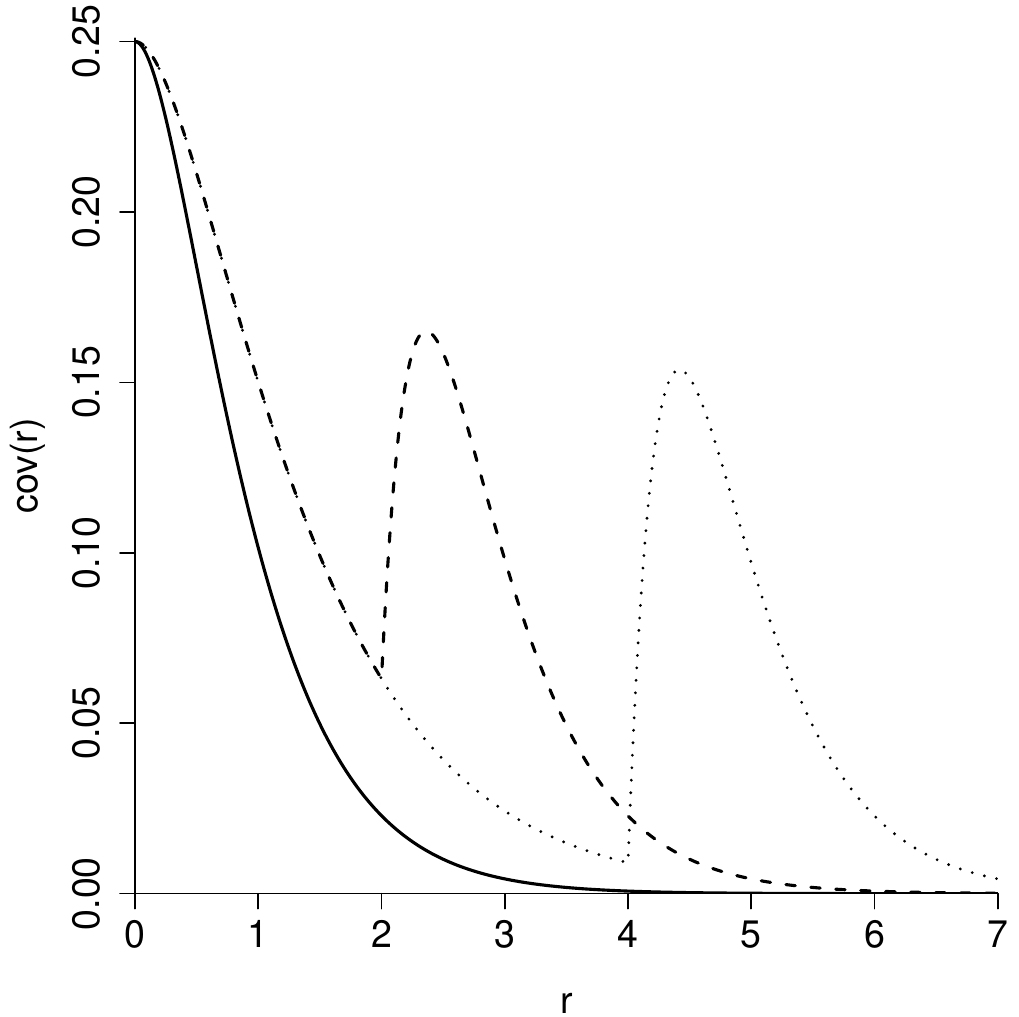}
\quad
\includegraphics[height=6cm]{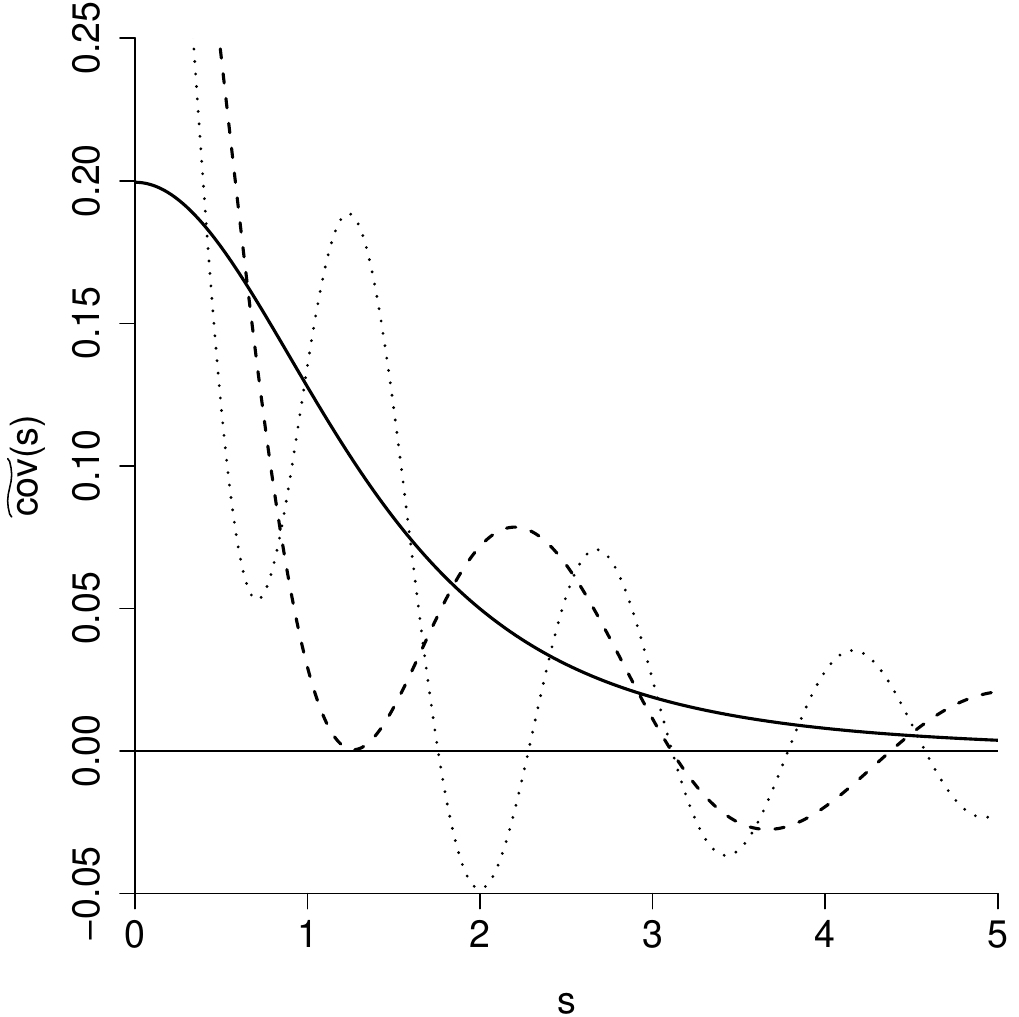}
\caption{Mark covariance function $\mathrm{cov}$ (left) and its Fourier transform $\widetilde{\mathrm{cov}}$ (right) for the model in Example \ref{ex:covEps} for $\lambda=1$ and $\varepsilon=0$ (solid), $\varepsilon=1$ (dashed), and $\varepsilon=2$ (dotted).}
\label{fig:covEps}
\end{center}
\end{figure}
\end{example}
Our major example is analysed within an own section since some results might be of interest not only to the field of random marked sets but also to the theory of positive definite functions.
%
%
\section{Gaussian random fields exceeding $t\in\mathbb{R}$}\label{sec:Gaussian}
Let $Z$ be a stationary and isotropic centered unit variance Gaussian random field in $\mathbb{R}^d$.  Then, for $t\in\mathbb{R}$, we define
\begin{align*}
\Xi_t&=\{x\in\mathbb{R}^d:\,Z(x)\geq t\}.
\end{align*}
If, in particular, $Z$ is almost surely continuous \cite{ref:Adler1981,ref:AdlerTaylor2007} then $\Xi_t$ is almost surely closed, i.\,e., $(\Xi_t,Z)$ is a random marked closed set. Note that $\Xi_t$ is a so-called excursion set which has been extensively studied in the literature, see \cite{ref:Adler1981,ref:AdlerTaylor2007} and the references therein. 

Since $Z$ is assumed to be both stationary and isotropic, its covariance
function,  
$
\mathrm{Cov}(x,y)=\mathbb{E}[Z(x)Z(y)]$, $x,y\in\mathbb{R}^d$,
is translation and rotation invariant, i.\,e., there exists a function
$R:[0,\infty)\rightarrow\mathbb{R}$ such that 
$\mathrm{Cov}(x,y)=R(\|x-y\|)$ for $x,y\in\mathbb{R}^d$.

First, we consider the case $t=0$.
%
\begin{theorem}\label{thm:Gauss_t=0}
Let $Z$ be a stationary and isotropic centered unit variance Gaussian random field in $\mathbb{R}^d$ with covariance function given by $R:[0,\infty)\rightarrow\mathbb{R}$.
Then, for $r\in[0,\infty)$, the second-order characteristics of $(\Xi_0,Z)$ are given by
\begin{align}
E(r)&=\sqrt{\frac{\pi}{2}}\,\frac{1+R(r)}{\arcsin(R(r))+\frac{\pi}{2}},\label{eqn:Gauss_t=0:E}\\
\Cov(r)&=R(r)+\frac{\sqrt{1-R(r)^2}}{\arcsin(R(r))+\frac{\pi}{2}}-\frac{\pi}{2}\frac{\left(1+R(r)\right)^2}{\left(\arcsin(R(r))+\frac{\pi}{2}\right)^2},\\
\gamma(r)&=(1-R(r))\left(1-\frac{\sqrt{1-R(r)^2}}{\arcsin(R(r))+\frac{\pi}{2}}\right),\\
k_{mm}(r)&=\frac{\pi}{2}\left(R(r)+\frac{\sqrt{1-R(r)^2}}{\arcsin(R(r))+\frac{\pi}{2}}\right),\\
\Cor(r)&=\tfrac{R(r)(\arcsin(R(r))+\frac{\pi}{2})^2+\sqrt{1-R(r)^2}(\arcsin(R(r))+\frac{\pi}{2})-\frac{\pi}{2}(1+R(r))^2}{(\arcsin(R(r))+\frac{\pi}{2})^2+R(r)\sqrt{1-R(r)^2}(\arcsin(R(r))+\frac{\pi}{2})-\frac{\pi}{2}(1+R(r))^2}.
\end{align}
\end{theorem} 
Obviously, each of the second-order characteristics of $(\Xi_0,Z)$ is a
continuous transform of $R$. In particular, this means that continuity of $R$
is preserved. Vice versa, due to the monotonicity of the transform for $\Cov$
(see Theorem \ref{thm:f_0} below) already from $\Cov$ it can be deduced
whether or not $R$ is continuous at the origin \cite{ref:Cressie1993}. 

Since, for every stationary random-field model, $E(h)$, $h\in\RR^d$, is
constant, equation (\ref{eqn:Gauss_t=0:E})  implies that
$(\Xi_0,Z)$ is not a random-field model, i.\,e., there does not exist a random
field in $\mathbb{R}^d$ whose second-order characteristics coincide with that
of $(\Xi_0,Z)$ unless $R$ is constant. It is therefore quite surprising to see
that we are not able 
to falsify that $(\Xi_0,Z)$ is a random-field model by using the mark
covariance function or the mark correlation function of $(\Xi_0,Z)$.
%
\begin{theorem}\label{thm:f_0}
The functions $f_0:[-1,1]\rightarrow\mathbb{R}$,
\begin{align*}
f_0(\rho)&=\rho+\frac{\sqrt{1-\rho^2}}{\arcsin \rho+\tfrac{\pi}{2}}-\tfrac{\pi}{2}\frac{(1+\rho)^2}{(\arcsin \rho+\tfrac{\pi}{2})^2},
\end{align*}
and $g_0:[-1,1]\rightarrow\mathbb{R}$,
\begin{align*}
g_0(\rho)&=\frac{\rho(\arcsin \rho+\tfrac{\pi}{2})^2+\sqrt{1-\rho^2}(\arcsin \rho+\tfrac{\pi}{2})-\tfrac{\pi}{2}(1+\rho)^2}{(\arcsin \rho+\tfrac{\pi}{2})^2+\rho\sqrt{1-\rho^2}(\arcsin \rho+\tfrac{\pi}{2})-\tfrac{\pi}{2}(1+\rho)^2},
\end{align*}
are  absolutely monotone on $[0,1]$, i.\,e., they have
only nonnegative derivatives there.
\end{theorem}
%
\begin{corollary}\label{cor:covcor_0PosDef}
$\mathrm{cov}(r)$ and $\mathrm{cor}(r)$ are positive definite functions.
\end{corollary}
However, the following proposition shows that the $k_{mm}$-function of
$(\Xi_0,Z)$ is not positive definite, in general.
%
\begin{proposition}\label{prop:kmm0}
Let $Z$ be a stationary and isotropic centered unit variance Gaussian random field in $\mathbb{R}^d$ with covariance function given by a continuous function $R:[0,\infty)\rightarrow\mathbb{R}$. Then the $k_{mm}$-function of $(\Xi_0,Z)$ is positive definite if and only if $R\equiv1$.
\end{proposition}
\begin{proof}
Let $q(\rho)=\rho+\sqrt{1-\rho^2}\,(\arcsin(\rho)+\pi/2)^{-1}$.
Then $q(1)=1$ and
\begin{align*}
q(\rho)&\geq \rho+\pi^{-1}\sqrt{1-\rho^2}>1,\quad \rho\in((\pi^2-1)/(\pi^2+1),1).
\end{align*}
Hence, $k_{mm}(r)$ is not a positive definite function if $R\not\equiv1$ \cite[Theorem 1.4.1]{ref:Sasvari1994}.
\end{proof}
The mark variogram of $(\Xi_0,Z)$ is in general not conditionally negative definite, which can be seen as follows. Consider $(\Xi_0,Z)$ for dimension $d=1$ and $R(r)=\cos(r)$. Since $\gamma$ is conditionally negative definite if and only if $\mathrm{e}^{-s\gamma}$ is positive definite for all $s>0$ \cite[Theorem 6.1.9]{ref:Sasvari1994} it suffices to show that $\mathrm{e}^{-\gamma}$ is not positive definite. $\mathrm{e}^{-\gamma(r)}$ inherits $2\pi$-periodicity from $\cos(r)$, and hence, it is positive definite if and only if its Fourier coefficients are nonnegative. Numerical calculations yield that the first Fourier coefficient is nearby -0.03364.
\bigskip

\begin{figure}
\begin{center}
\includegraphics[height=6cm]{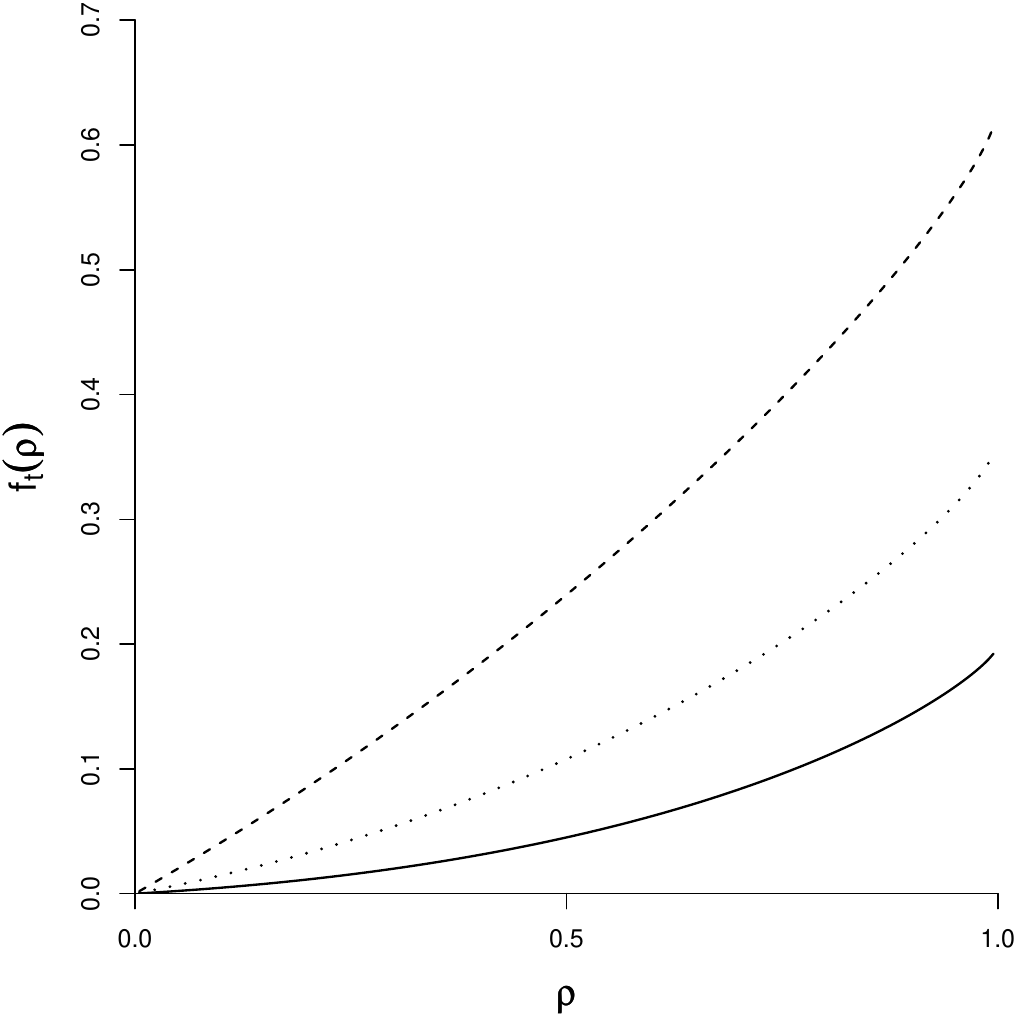}
\quad
\includegraphics[height=6cm]{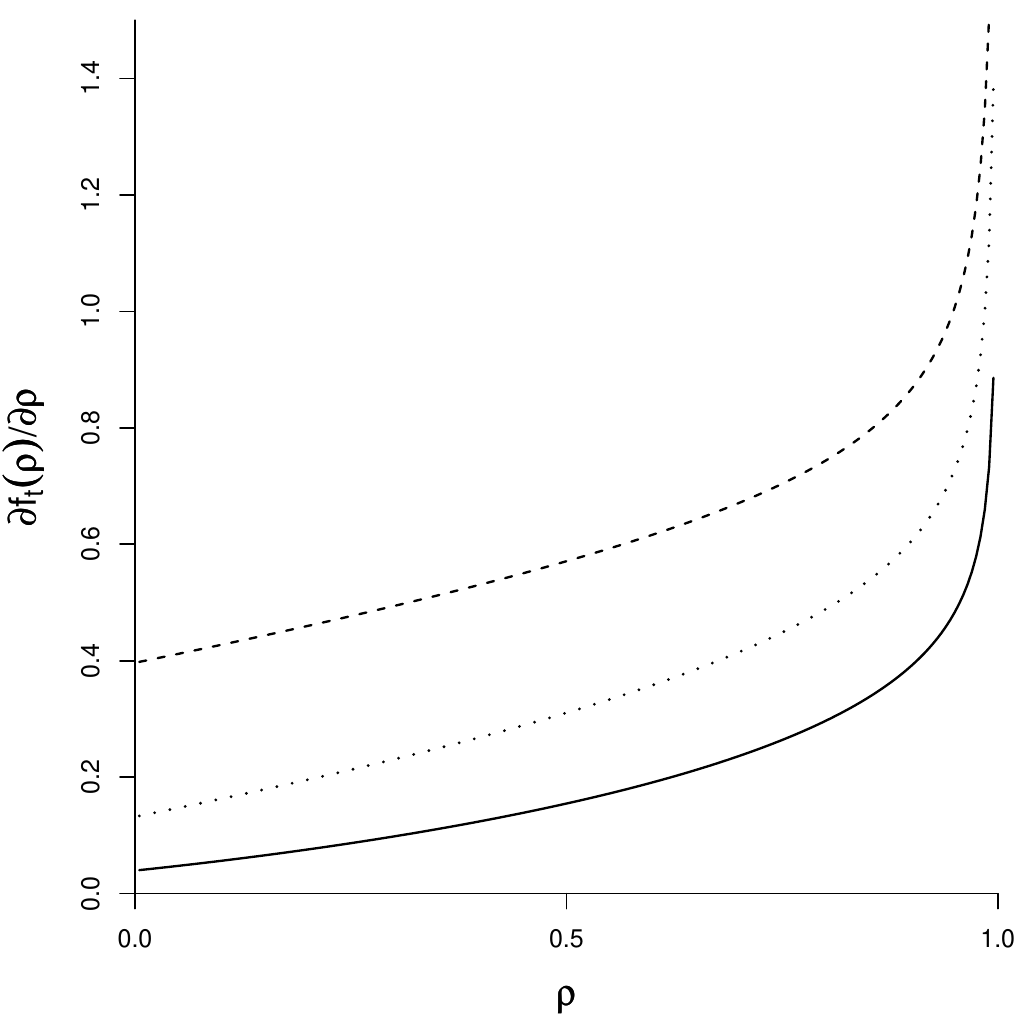}
\caption{$f_t$ (left) and $f_t'$ (right) for $t=-1$ (dashed), $t=0$ (dotted) and $t=1$ (solid).}
\label{fig:covt}
\end{center}
\end{figure}

Now we switch over to the more general case $t\in\mathbb{R}$. Unfortunately,
unlike the case $t=0$, we cannot express all the second-order characteristics
of $(\Xi_t,Z)$ in closed form. In particular, for a stationary and isotropic
centered unit variance Gaussian random field $Z$ in $\mathbb{R}^d$ with
covariance function given by $R:[0,\infty)\rightarrow\mathbb{R}$, we have 
\begin{align*}
\mathbb{P}(o,h\in\Xi_t)&=
\int\limits_0^{R(\|h\|)}\varphi(t,t,s)\,\mathrm{d}s+\Psi(t)^2,
\end{align*} 
see \cite[Eqn. (10.8.3)]{ref:CramerLeadbetter1967}. Here, 
\[\Psi(t)=\int\limits_t^{\infty} \varphi(s) \,\mathrm{d}s,
\qquad
\varphi(t) =
\frac{1}{\sqrt{2\pi}} \mathrm{e}^{-\frac{t^2}{2}}, \qquad t \in\RR,
\] 
denotes the tail probability function of the standard Gaussian distribution.
 By $\varphi(x,y,\rho)$ we denote the density of the bivariate Gaussian
 distribution with unit variances and correlation $\rho$.
 In the following we concentrate on the mark covariance function and the mark variogram of $(\Xi_t,Z)$ and write
\begin{align*}
P_t(\rho)&=\int_0^{\rho} \varphi(t, t, s)\,\mathrm{d}s+\Psi(t)^2, \quad \rho\in[-1,1].
\end{align*}
%
\begin{lemma}\label{lem:E_t+C_t}
Let $Z$ be a stationary and isotropic centered unit variance Gaussian random
field in $\mathbb{R}^d$ with covariance function given by
$R:[0,\infty)\rightarrow\mathbb{R}$. 
Then, for $t\in\mathbb{R}$ and $h\in\mathbb{R}^d$, we have
\begin{align*}
\mathbb{E}[Z(o)\mathbf{1}_{\Xi_t}(o)\mathbf{1}_{\Xi_t}(h)]
&\,= E_t(R(\|h\|)),
\quad E_t(\rho) = 
\varphi(t)\big(\rho+1\big)\Psi\left(t\sqrt{\frac{1-\rho}{1+\rho}}\right),
\end{align*}
\begin{align*}
&\mathbb{E}[Z(o)Z(h)\mathbf{1}_{\Xi_t}(o)\mathbf{1}_{\Xi_t}(h)]
= C_t(R(\|h\|))
\end{align*}
where
\begin{align*}
C_t(\rho)
&\;=
(1-\rho^2)\varphi(t,t,s)+2\rho t\varphi(t)\Psi\left(t\sqrt{(1-\rho)/(1+\rho)}\right)+\rho P_t(\rho),
\end{align*}
and
\begin{align*}
&\mathbb{E}[Z(o)^2\mathbf{1}_{\Xi_t}(o)\mathbf{1}_{\Xi_t}(h)]
= V_t(R(\|h\|))
\end{align*}
where
\begin{align*}
V_t(\rho)
&\;=
\rho(1-\rho^2)\varphi(t,t,s)+(1+\rho^2) t\varphi(t)\Psi\left(t\sqrt{(1-\rho)/(1+\rho)}\right)+P_t(\rho).
\end{align*}
\end{lemma}
Then the mark covariance function of $(\Xi_t,Z)$ is given by
\begin{equation}\label{eqn:f_t}
\mathrm{cov}(r)=f_t(R(r)),\quad f_t(\rho)=\frac{C_t(\rho)}{P_t(\rho)}-\frac{E_t(\rho)^2}{P_t(\rho)^2},
\end{equation}
and the mark variogram of $(\Xi_t,Z)$ is given by
\begin{equation}\label{eqn:v_t}
\gamma(r)=v_t(R(r)),\quad v_t(\rho)=\frac{V_t(\rho)-C_t(\rho)}{P_t(\rho)}.
\end{equation}
There is strong evidence that also the mark covariance function of $(\Xi_t,Z)$, $t\neq0$, is positive definite for a certain class of Gaussian random fields $Z$. Figure \ref{fig:covt} shows $f_t(\rho)$ and $f_t'(\rho)$ for several $t$, indicating that for these $t$ the functions $f_t(\rho)$ are both increasing and convex for $\rho\in[0,1]$. Hence, if this is really true, for instance P{\'o}lya's criterion \cite{ref:Polya1949} would imply that, for any continuous and convex function $R:[0,\infty)\rightarrow\mathbb{R}$ satisfying $R(0)=1$ and $\lim_{r\rightarrow\infty}R(r)=0$, the function $f_t(R(|\cdot|))$ is positive definite on $\mathbb{R}$.
%
%
%
\section{Differentiability at 0}\label{sec:Differentiability}
In this section we continue the example of the preceding section and show that the corresponding mark covariance function has a right-hand derivative at 0 which does not vanish. 
%
\begin{lemma}\label{lem:covExcursionSet}
Let $Z$ be a stationary and isotropic centered unit variance Gaussian random
field in $\mathbb{R}^d$ with continuous covariance function given by
$R:[0,\infty)\rightarrow\mathbb{R}$ excluding $R\equiv1$. Let
$C_{\Xi_t}(\|h\|)=\mathbb{P}(o,h\in\Xi_t)$
be the set covariance of the excursion
set $\Xi_t$. If $Z$ is mean-square differentiable then 
\begin{align*}
C_{\Xi_t}'(0+)&=-\frac{\varphi(t)}{\sqrt{2\pi}}\sqrt{-R''(0+)} < 0,
\end{align*}
otherwise $C_{\Xi_t}'(0+)=-\infty$.
\end{lemma}
%
\begin{theorem}\label{thm:covtDiff}
Under the assumptions of Lemma \ref{lem:covExcursionSet}
the mark covariance function $\mathrm{cov}(r)=f_t(R(r))$ of $(\Xi_t,Z)$ has a negative right-hand derivative at $r=0$, in particular, we have
\begin{align*}
\mathrm{cov}'(0+)&=-\frac{(t^2-1)\Psi(t)^2-3t\varphi(t)\Psi(t)+2\varphi(t)^2}{\Psi(t)^3}\cdot\frac{\varphi(t)\sqrt{-R''(0+)}}{\sqrt{2\pi}}<0
\end{align*}
in case $Z$ is mean-square differentiable, and $\mathrm{cov}'(0+)=-\infty$ otherwise. In either case, the right-hand derivative of the mark variogram $\gamma(r)=v_t(R(r))$ of $(\Xi_t,Z)$ at $r=0$ is $-R'(0+)$.
\end{theorem}


\section{Proofs}\label{sec:Proofs}
%
First, we prove  Lemma \ref{lem:E_t+C_t} of Section \ref{sec:Gaussian}
which is needed in the proof of Theorem \ref{thm:Gauss_t=0}.
\subsection{Proof of Lemma \ref{lem:E_t+C_t}}\label{subsec:proof:E_t+C_t}
We will repeatedly apply the identities
\begin{equation}\label{eqn:id1}
\varphi'(y)=-y\varphi(y),
\end{equation}
\begin{equation}\label{eqn:id2}
\varphi(y,x,\rho)=\varphi(x,y,\rho)
=\frac{1}{\sqrt{1-\rho^2}}\varphi(x)\varphi\bigg(\frac{y-\rho x}{\sqrt{1-\rho^2}}\bigg),
\end{equation}
\begin{equation}\label{eqn:id3}
\int_t^{\infty}\frac{1}{\sqrt{1-\rho^2}}\varphi\bigg(\frac{t-\rho x}{\sqrt{1-\rho^2}}\bigg)\varphi(x)\,\dd x
=\varphi(t)\Psi\bigg(t\sqrt{\frac{1-\rho}{1+\rho}}\bigg),
\end{equation}
and
\begin{align}
P_t(\rho)&=\int_t^{\infty}\int_t^{\infty}\!\!\frac{1}{\sqrt{1-\rho^2}}\varphi\bigg(\frac{x_1-\rho x_2}{\sqrt{1-\rho^2}}\bigg)\varphi(x_2)\,\dd x_1\,\dd x_2\nonumber\\
&=\int_t^{\infty}\!\!\Psi\bigg(\frac{t-\rho x_2}{\sqrt{1-\rho^2}}\bigg)\varphi(x_2)\,\dd x_2.\label{eqn:id4}
\end{align}
By (\ref{eqn:id2}), the change of variables $x_1=\sqrt{1-\rho^2}\,y+\rho x_2$, (\ref{eqn:id1}) and (\ref{eqn:id3}) we have
\begin{align*}
E_t(\rho)&=\int_t^{\infty}\int_t^{\infty}x_1\frac{1}{\sqrt{1-\rho^2}}\,\varphi\bigg(\frac{x_1-\rho x_2}{\sqrt{1-\rho^2}}\bigg)\varphi(x_2)\,\dd x_1\,\dd x_2\\
&=\int_t^{\infty}\int_{\frac{t(1-\rho)}{\sqrt{1-\rho^2}}}^{\infty}(\sqrt{1-\rho^2}\,y+\rho x_2)\varphi(y)\varphi(x_2)\,\dd y\,\dd x_2\\
&=\int_t^{\infty}\sqrt{1-\rho^2}\,\varphi\bigg(\frac{t-\rho x_2}{\sqrt{1-\rho^2}}\bigg)\varphi(x_2)+\rho x_2\Psi\bigg(\frac{t-\rho x_2}{\sqrt{1-\rho^2}}\bigg)\varphi(x_2)\,\dd x_2\\
&=(1-\rho^2)\varphi(t)\Psi\bigg(t\sqrt{\frac{1-\rho}{1+\rho}}\bigg)+\int_t^{\infty}\rho x_2\Psi\bigg(\frac{t-\rho x_2}{\sqrt{1-\rho^2}}\bigg)\varphi(x_2)\,\dd x_2.
\end{align*}
Integration by parts, (\ref{eqn:id1}) and (\ref{eqn:id3}) finally yield
\begin{align*}
&\int_t^{\infty}\rho x_2\Psi\bigg(\frac{t-\rho x_2}{\sqrt{1-\rho^2}}\bigg)\varphi(x_2)\,\dd x_2\\
&\quad=\rho\varphi(t)\Psi\bigg(t\sqrt{\frac{1-\rho}{1+\rho}}\bigg)+\int_t^{\infty}\frac{\rho^2}{\sqrt{1-\rho^2}}\varphi\bigg(\frac{t-\rho x_2}{\sqrt{1-\rho^2}}\bigg)\varphi(x_2)\,\dd x_2\\
&\quad=(\rho+\rho^2)\varphi(t)\Psi\bigg(t\sqrt{\frac{1-\rho}{1+\rho}}\bigg).
\end{align*}
Similarly to the calculations for $E_t(\rho)$ we have
\begin{align*}
C_t(\rho)&=\int_t^{\infty}\sqrt{1-\rho^2}\,x_2\varphi\bigg(\frac{t-\rho x_2}{\sqrt{1-\rho^2}}\bigg)\varphi(x_2)\,\dd x_2\\
&\qquad+\int_t^{\infty}\rho x_2^2\Psi\bigg(\frac{t-\rho x_2}{\sqrt{1-\rho^2}}\bigg)\varphi(x_2)\,\dd x_2.
\end{align*}
Integration by parts, (\ref{eqn:id1}) and (\ref{eqn:id4}) yield for the second integral
\begin{align*}
&\int_t^{\infty}\rho x_2^2\Psi\bigg(\frac{t-\rho x_2}{\sqrt{1-\rho^2}}\bigg)\varphi(x_2)\,\dd x_2\\
&\quad=\rho t\varphi(t)\Psi\bigg(t\sqrt{\tfrac{1-\rho}{1+\rho}}\bigg)+\rho P_t(\rho)+\int_t^{\infty}\frac{\rho^2}{\sqrt{1-\rho^2}}\,x_2\varphi\bigg(\tfrac{t-\rho x_2}{\sqrt{1-\rho^2}}\bigg)\varphi(x_2)\,\dd x_2.
\end{align*}
Hence, we have
\begin{equation*}
C_t(\rho)=\int_t^{\infty}\frac{x_2}{\sqrt{1-\rho^2}}\,\varphi\bigg(\tfrac{t-\rho x_2}{\sqrt{1-\rho^2}}\bigg)\varphi(x_2)\,\dd x_2+\rho t\varphi(t)\Psi\bigg(t\sqrt{\tfrac{1-\rho}{1+\rho}}\bigg)+\rho P_t(\rho).
\end{equation*}
By (\ref{eqn:id2}), the change of variables $x_2=\sqrt{1-\rho^2}\,y+\rho t$, (\ref{eqn:id1}) and again (\ref{eqn:id2}) we obtain for the remaining integral
\begin{align*}
&\int_t^{\infty}\frac{x_2}{\sqrt{1-\rho^2}}\,\varphi\bigg(\frac{t-\rho x_2}{\sqrt{1-\rho^2}}\bigg)\varphi(x_2)\,\dd x_2\\
&\quad=\varphi(t)\int_t^{\infty}\frac{x_2}{\sqrt{1-\rho^2}}\,\varphi\bigg(\frac{x_2-\rho t}{\sqrt{1-\rho^2}}\bigg)\,\dd x_2\\
&\quad=\varphi(t)\int_{\frac{t(1-\rho)}{\sqrt{1-\rho^2}}}^{\infty}(\sqrt{1-\rho^2}\,y+\rho t)\,\varphi(y)\,\dd y\\
&\quad=\sqrt{1-\rho^2}\,\varphi(t)\varphi\bigg(\frac{t-\rho t}{\sqrt{1-\rho^2}}\bigg)+\rho t\varphi(t)\Psi\bigg(t\sqrt{\frac{1-\rho}{1+\rho}}\bigg)\\
&\quad=(1-\rho^2)\varphi(t,t\rho)+\rho t\varphi(t)\Psi\bigg(t\sqrt{\frac{1-\rho}{1+\rho}}\bigg).
\end{align*}
Finally, using the identities (\ref{eqn:id1})--(\ref{eqn:id4}), the result for $V_t(\rho)$ can be verified similarly.
\subsection{Proof of Theorem \ref{thm:Gauss_t=0}}
As $\varphi(0,0,\rho)=(2\pi)^{-1}(1-\rho^2)^{-1/2}$ and
\begin{align*}
\mathbb{P}(Z(o)\geq0,\,Z(r)\geq0)&=\frac{1}{2\pi}\left(\arcsin(R(r))+\frac{\pi}{2}\right)
\end{align*}
the formulae for $E$, $\Cov$, $\gamma$ and $\Cor$ follow immediately from Lemma
\ref{lem:E_t+C_t}. Finally, $\overline{m}=\mathbb{E}[Z(o)\mid Z(o)\geq0]=\sqrt{2/\pi}$ yields the result for $k_{mm}$.
%
\subsection{Proof of Theorem \ref{thm:f_0}}
The idea of the proof is to show that the coefficients of the Taylor expansion
of $f_0$ and $g_0$ are all nonnegative. We give the proof for $f_0$ only;
for $g_0$, the same techniques can be applied.

Since $\arcsin$ can be extended to the complex plain $\mathbb{C}$ by
\begin{align*}
\arcsin z&=-\mathrm{i}\ln(\mathrm{i}z+\sqrt{1-z^2}),\quad z\in\mathbb{C},
\end{align*}
with branch cuts $(-\infty,-1)$ and $(1,\infty)$,
 also $f_0$ can be extended to $\mathbb{C}$. 
We restrict our attention to 
\begin{align*}
f(z)&=\frac{\sqrt{1-z^2}}{\arcsin z+\tfrac{\pi}{2}}-\tfrac{\pi}{2}\frac{(1+z)^2}{(\arcsin z+\tfrac{\pi}{2})^2}.
\end{align*}
Since $\arcsin z=-\pi/2$ is satisfied only for $z=-1$ and the denumerators of all derivatives of $f$ consist of a product of integer powers  
of $\arcsin z+\pi/2$ and $\sqrt{1-z^2}$, $f$ is not only holomorphic in the open unit disc $D=\{z\in\mathbb{C}:|z|<1\}$
but also on $\mathbb{C}\setminus((-\infty,-1]\cup[1,\infty))$. This holds because $\sqrt{z}$ is holomorphic everywhere except on $(-\infty,0]$ where $\sqrt{z}$ is not even continuous. The expansion of $\arcsin z$ around $z=-1$
\begin{align*}
\arcsin z&=-\tfrac{\pi}{2}+\sqrt{2}\sqrt{1+z}\big(1+O(1+z)\big),\quad|z|<1,
\end{align*}
shows that the derivatives of any order of $f$ can be continuously extended in $z=-1$ as $z\rightarrow-1$, $|z|<1$. Hence, the only branching point of $f$ in $\overline{D}$ is $z=1$. From the expansion of $f$ around $z=1$
\begin{align*}
f(z)&=-\tfrac{2}{\pi}+c_1(1-z)^{\tfrac{1}{2}}+\left(\tfrac{2}{\pi}+\tfrac{2}{\pi^2}-\tfrac{12}{\pi^3}\right)(1-z)
+c_2(1-z)^{\tfrac{3}{2}}\\
&\quad+O((1-z)^2),\quad|z|<1,
\end{align*}
with
\begin{align*} \nonumber
&\qquad c_1=-\sqrt{2}\left(\tfrac{4}{\pi^2}-\tfrac{1}{\pi}\right),\quad c_2=\sqrt{2}\left(\tfrac{-1}{4\pi}+\tfrac{11}{3\pi^2}+\tfrac{2}{\pi^3}-\tfrac{16}{\pi^4}\right),
\end{align*}
we see that the asymptotic behaviour of $f$ around $z=1$ is essentially described by
\begin{align*}
g(z)&=c_1(1-z)^{\tfrac{1}{2}}+c_2(1-z)^{\tfrac{3}{2}}.
\end{align*}
At first, we derive an upper bound for the coefficients $b_n$ of the Taylor expansion of 
\begin{align*}
h(z)&=f(z)-g(z)
\end{align*}
at $z=0$. The function 
$h$ is holomorphic in $D$ and by construction twice continuously differentiable in $\overline{D}$ since $h^{(2)}(z)$ behaves around $z=1$ like
$c+O((1-z)^{\tfrac{1}{2}})$, $|z|<1$. As a consequence of Cauchy's integral formula we have
\begin{align*}
b_n&=\frac{h^{(n)}(0)}{n!}=\frac{1}{2\pi\mathrm{i}}\int\limits_{\partial D_{1-\varepsilon}}\frac{h(\zeta)}{\zeta^{n+1}}\,\mathrm{d}\zeta,\quad n=0,1,2,\ldots,
\end{align*}
for all $\varepsilon\in(0,1)$ and $D_{1-\varepsilon}=\{z\in\mathbb{C}:|z|<1-\varepsilon\}$. After twice partially integrating we have
\begin{align*}
b_n&=\frac{1}{2\pi\mathrm{i}}\frac{1}{n(n-1)}\int\limits_{\partial D_{1-\varepsilon}}\frac{h^{(2)}(\zeta)}{\zeta^{n-1}}\,\mathrm{d}\zeta,\quad n=2,3,4,\ldots.
\end{align*}
Hence, for $n\geq2$, we obtain 
\begin{align*}
|b_n|&\leq\frac{1}{2\pi}\frac{1}{n(n-1)}(1-\varepsilon)^{1-n}\max_{\varphi\in[0,2\pi]}\left|h^{(2)}((1-\varepsilon)\mathrm{e}^{\mathrm{i}\varphi})\right|\cdot2\pi(1-\varepsilon).
\end{align*}
Since $h^{(2)}(z)$ is continuous in $\overline{D}$ we can conclude that for $n\geq2$ 
\begin{align*}
|b_n|&\leq\frac{1}{n(n-1)}\max_{\varphi\in[0,2\pi]}\left|h^{(2)}(\mathrm{e}^{\mathrm{i}\varphi})\right|<\frac{1}{n(n-1)}0.182
\end{align*}
by Lemma \ref{lem:Estimate} below.
In the following we derive a lower bound for the coefficients $a_n$ of the
Taylor expansion of $g(z)$ at $z=0$. 
Note that we have the series expansions 
\begin{align*}
(1-z)^{\tfrac{1}{2}}&=\sum_{n=0}^{\infty}(-1)^n\binom{\frac{1}{2}}{n}z^n=1-\sum_{n=1}^{\infty}\frac{(2n-2)!}{n!(n-1)!2^{2n-1}}z^n,\quad|z|<1,
\end{align*}
and
\begin{align*}
(1-z)^{\tfrac{3}{2}}&=1-\tfrac{3}{2}z+3\sum_{n=2}^{\infty}\frac{(2n-4)!}{n!(n-2)!2^{2n-2}}z^n,\quad|z|<1.
\end{align*}
Hence, for $n\geq2$, we have
\begin{align*}
a_n&=-c_1\cdot\frac{(2n-2)!}{n!(n-1)!2^{2n-1}}+c_2\cdot3\,\frac{(2n-4)!}{n!(n-2)!2^{2n-2}}.
\end{align*}
Since $c_1<0$ and $c_2>0$, we obtain
\begin{align*}
a_n&>-c_1\cdot\frac{(2n-2)!}{n!(n-1)!2^{2n-1}}
\ge 
-c_1
\frac{\mathrm{e}}{2\sqrt{\pi}}\left(\frac{n-1}{n}\right)^{n-1}n^{-\frac{3}{2}}\mathrm{e}^{-\frac{1}{12n}-\frac{1}{12(n-1)}}
\end{align*}
by Stirling's approximation.
At least for $n\geq6$ we have
\begin{align*}
e \left(\tfrac{n-1}{n}\right)^n \cdot
\left(1+\tfrac{1}{n-1}\right)\mathrm{e}^{-\frac{1}{12n}-\frac{1}{12(n-1)}}&\geq\left(1+\tfrac{1}{n-1}\right)\mathrm{e}^{-\frac{1}{6(n-1)}}
\geq1.
\end{align*}
Thus,
\begin{align*}
&a_n > \frac{4-\pi}{\pi^2\sqrt{2\pi}}\,n^{-\frac{3}{2}},
\qquad n\ge 6.
\end{align*}
Combining this with the upper bound for $|b_n|$ it is easily seen that
 the $n$th coefficient of the
Taylor expansion of $f_0$ at $x=0$ is nonnegative for all $n\geq30$. The
remaining coefficients of order $n=0,1,\ldots,29$ can be determined from a
series expansion and are easily evaluated to be nonnegative.

The same techniques can be applied to prove the assertion for $g_0$.
\begin{lemma}\label{lem:Estimate}
We have
\begin{align*}
\max_{\varphi\in[0,2\pi]}\left|h^{(2)}(\mathrm{e}^{\mathrm{i}\varphi})\right|&<0.182.
\end{align*}
\end{lemma}
\begin{proof}
Figure \ref{fig:Proof_abs_mon} indicates that this maximum is achieved for $\varphi=0$. From the expansion of $f$ around $z=1$ we obtain the value of this maximum: It is twice the absolute value of the coefficent of $(1-z)^2$, i.\,e., we have
\begin{align*}
\max_{\varphi\in[0,2\pi]}\left|h^{(2)}(\mathrm{e}^{\mathrm{i}\varphi})\right|&=|h^{(2)}(1)|=\left(\tfrac{-1}{\pi}-\tfrac{2}{3\pi^2}+\tfrac{20}{\pi^3}+\tfrac{8}{\pi^4}-\tfrac{80}{\pi^5}\right)<0.08.
\end{align*}
Formally, we can actually prove that an upper bound is 0.182. Since the conjugate-complex of $h^{(2)}(\mathrm{e}^{\mathrm{i}\varphi})$
equals $h^{(2)}(\mathrm{e}^{-\mathrm{i}\varphi})$ we can restrict to the
interval $[0,\pi]$. 

First we make an estimate for $[\tfrac{\pi}{3},\pi]$ then for $[0,\tfrac{\pi}{3}]$. From the series expansion
\begin{align}\label{eqn:sqrt1-z}
(1-z)^{\tfrac{1}{2}}&=1-\sum_{n=1}^{\infty}\frac{(2n-2)!}{n!(n-1)!2^{2n-1}}z^n,\quad|z|<1,
\end{align}
we obtain the two series expansions at $z=-1$
\begin{align*}
(1-z)^{\tfrac{1}{2}}&=
\sqrt{2}\sqrt{1-\frac{1+z}{2}}=\sqrt{2}\left(1-\sum_{n=1}^{\infty}\frac{(2n-2)!}{n!(n-1)!2^{2n-1}}\left(\frac{1+z}{2}\right)^n\right)
\\&=\sum_{n=0}^{\infty}v_n(1+z)^n\hfill
\end{align*}
and
\begin{align*}
\arcsin
z+\tfrac{\pi}{2}&=\sqrt{2}\sum_{n=0}^{\infty}\frac{(2n)!}{n!n!2^{3n+1}}\,\frac{1}{n+\frac{1}{2}}(1+z)^{n+\frac{1}{2}}
\\&=(1+z)^{\frac{1}{2}}\sum_{n=0}^{\infty}u^{(1)}_n(1+z)^{n},
\end{align*}
which are convergent for $|1+z|<2$. From the latter we also obtain according series expansions
\begin{align*}
(\arcsin z+\tfrac{\pi}{2})^k&=(1+z)^{\frac{k}{2}}\sum_{n=0}^{\infty}u^{(k)}_n(1+z)^{n},\quad k=2,3,4.
\end{align*}
Let
\begin{align*}
P_N^{(k)}(z)&=(1+z)^{\frac{k}{2}}\sum_{n=0}^{N}u^{(k)}_n(1+z)^{n},\qquad k=1,2,3,4
\\
p_N^{(k)}(z)&=\frac{\arcsin
  z+\tfrac{\pi}{2}-P_N^{(k)}(z)}{(1+z)^{\frac{7}{2}}},
\qquad k=1,2,3,4
,\end{align*}
and
\begin{align*}
 Q_N(z)=\sum_{n=0}^{N}v_n(1+z)^{n},
\qquad q_N(z)=\frac{(1-z)^{\tfrac{1}{2}}-Q_N(z)}{(1+z)^{\frac{7}{2}}}
.\end{align*}
We write $h^{(2)}(z)$ as a fraction with numerator
\begin{align*}
h_{\mathrm{num}}(z)&=2(\arcsin z+\tfrac{\pi}{2})(1+z)(1-z)-(\arcsin z+\tfrac{\pi}{2})^3(1+z)(1-z)\\ \nonumber
&\quad-\pi(\arcsin z+\tfrac{\pi}{2})^2(1+z)^{\frac{3}{2}}(1-z)^{\frac{3}{2}}
\\&\quad  +z(\arcsin
z+\tfrac{\pi}{2})^2(1+z)^{\frac{1}{2}}(1-z)^{\frac{1}{2}}-z^2(\arcsin
z+\tfrac{\pi}{2})^3\\ \nonumber 
&\quad+\frac{c_1}{4}(\arcsin
z+\tfrac{\pi}{2})^4(1+z)^{\frac{3}{2}}-3\pi(1+z)^{\frac{5}{2}}(1-z)^{\frac{1}{2}}
\\&\quad+4\pi(\arcsin z+\tfrac{\pi}{2})(1+z)^2(1-z)\\ \nonumber
&\quad+\pi z(\arcsin z+\tfrac{\pi}{2})(1+z)^2-\frac{3\,c_2}{4}(\arcsin z+\tfrac{\pi}{2})^4(1+z)^{\frac{3}{2}}(1-z)
\end{align*}
and denumerator
\begin{align*}
h_{\mathrm{denum}}(z)&=(\arcsin z+\tfrac{\pi}{2})^4(1+z)^{\frac{3}{2}}(1-z)^{\frac{3}{2}}.
\end{align*}
By $h_{\mathrm{num}}^{\mathrm{approx}}(z)$ we denote the approximation of $h_{\mathrm{num}}(z)$ if we replace $(\arcsin z+\tfrac{\pi}{2})^k$ by $P_N^{(k)}(z)$ and $\sqrt{1-z}$ by $Q_N(z)$ within $h_{\mathrm{num}}(z)$. On the set $\{\mathrm{e}^{\mathrm{i}\varphi}:\varphi\in[\frac{\pi}{3},\pi]\}$ we can make the following estimate for the error
\begin{align*}
\eta(\varphi)&=\left|\frac{h_{\mathrm{num}}\left(\mathrm{e}^{\mathrm{i}\varphi}\right)-h_{\mathrm{num}}^{\mathrm{approx}}\left(\mathrm{e}^{\mathrm{i}\varphi}\right)}{h_{\mathrm{denum}}\left(\mathrm{e}^{\mathrm{i}\varphi}\right)}\right|.
\end{align*}
Since we have
\begin{align*}
|(\arcsin z+\tfrac{\pi}{2})^4|\geq4|1+z|^2
\end{align*}
we obtain
\begin{align*}
\eta(\varphi)&\leq\frac{1}{4}|1-\mathrm{e}^{\mathrm{i}\varphi}|^{-\frac{3}{2}}\bigg(
2|p_N^{(1)}(\mathrm{e}^{\mathrm{i}\varphi})|\cdot|1+\mathrm{e}^{\mathrm{i}\varphi}|\cdot|1-\mathrm{e}^{\mathrm{i}\varphi}|
\\&\quad
+|p_N^{(3)}(\mathrm{e}^{\mathrm{i}\varphi})|\cdot|1+\mathrm{e}^{\mathrm{i}\varphi}|\cdot|1-\mathrm{e}^{\mathrm{i}\varphi}|\\ \nonumber
&\quad+\pi|p_N^{(2)}(\mathrm{e}^{\mathrm{i}\varphi})|\cdot|1+\mathrm{e}^{\mathrm{i}\varphi}|^{\frac{3}{2}}\cdot|1-\mathrm{e}^{\mathrm{i}\varphi}|\cdot|Q_N(\mathrm{e}^{\mathrm{i}\varphi})|
\\&\quad
+\pi|P_N^{(2)}(\mathrm{e}^{\mathrm{i}\varphi})||1+\mathrm{e}^{\mathrm{i}\varphi}|^{\frac{3}{2}}|1-\mathrm{e}^{\mathrm{i}\varphi}|\cdot|q_N(\mathrm{e}^{\mathrm{i}\varphi})|\\ \nonumber
&\quad+\pi|p_N^{(2)}(\mathrm{e}^{\mathrm{i}\varphi})|\cdot|1+\mathrm{e}^{\mathrm{i}\varphi}|^{5}\cdot|1-\mathrm{e}^{\mathrm{i}\varphi}|\cdot|q_N(\mathrm{e}^{\mathrm{i}\varphi})|
\\&\quad
+|p_N^{(2)}(\mathrm{e}^{\mathrm{i}\varphi})|\cdot|1+\mathrm{e}^{\mathrm{i}\varphi}|^{\frac{1}{2}}\cdot|Q_N(\mathrm{e}^{\mathrm{i}\varphi})|\\ \nonumber
&\quad+|P_N^{(2)}(\mathrm{e}^{\mathrm{i}\varphi})|\cdot|1+\mathrm{e}^{\mathrm{i}\varphi}|^{\frac{1}{2}}\cdot|q_N(\mathrm{e}^{\mathrm{i}\varphi})|+|p_N^{(2)}(\mathrm{e}^{\mathrm{i}\varphi})|\cdot|1+\mathrm{e}^{\mathrm{i}\varphi}|^{4}\cdot|q_N(\mathrm{e}^{\mathrm{i}\varphi})|\\ \nonumber
&\quad+|p_N^{(3)}(\mathrm{e}^{\mathrm{i}\varphi})|+\frac{|c_1|}{4}|p_N^{(4)}(\mathrm{e}^{\mathrm{i}\varphi})|\cdot|1+\mathrm{e}^{\mathrm{i}\varphi}|^{\frac{3}{2}}+3\pi|1+\mathrm{e}^{\mathrm{i}\varphi}|^{\frac{5}{2}}\cdot|q_N(\mathrm{e}^{\mathrm{i}\varphi})|\\ \nonumber
&\quad+4\pi|1+\mathrm{e}^{\mathrm{i}\varphi}|^2\cdot|1-\mathrm{e}^{\mathrm{i}\varphi}|\cdot|p_N^{(1)}(\mathrm{e}^{\mathrm{i}\varphi})|+\pi|1+\mathrm{e}^{\mathrm{i}\varphi}|^2\cdot|p_N^{(1)}(\mathrm{e}^{\mathrm{i}\varphi})|\\ \nonumber
&\quad+\frac{3|c_2|}{4}|p_N^{(4)}(\mathrm{e}^{\mathrm{i}\varphi})|\cdot|1+\mathrm{e}^{\mathrm{i}\varphi}|^{\frac{3}{2}}\cdot|1-\mathrm{e}^{\mathrm{i}\varphi}|\bigg).
\end{align*}
Since $|1+\mathrm{e}^{\mathrm{i}\varphi}|$ is decreasing on $[0,\pi]$ we have 
\begin{align*}
|1+\mathrm{e}^{\mathrm{i}\varphi}|&\leq|1+\mathrm{e}^{\mathrm{i}\frac{\pi}{3}}|=\sqrt{3}
\end{align*}
on $[\frac{\pi}{3},\pi]$. The coefficients $u_n^{(k)}$, $k=1,2,3,4$, and $v_n$ are non-alternating for $n\geq1$. Hence, on $[\frac{\pi}{3},\pi]$, we have
\begin{align*}
|p_N^{(k)}(\mathrm{e}^{\mathrm{i}\varphi})|&\leq|p_N^{(k)}(\mathrm{e}^{\mathrm{i}\frac{\pi}{3}})|=:\delta_N^{(k)},\quad k=1,2,3,4,\\
|q_N(\mathrm{e}^{\mathrm{i}\varphi})|&\leq|q_N(\mathrm{e}^{\mathrm{i}\frac{\pi}{3}})|=:\varepsilon_N
\end{align*}
for all $N\geq1$. Furthermore, on $[\frac{\pi}{3},\pi]$, we have
\begin{align*}
|1-\mathrm{e}^{\mathrm{i}\varphi}|^{-1}&\leq1
\end{align*}
and, for $N\geq20$, we can estimate
\begin{align*}
|P_N^{(2)}(\mathrm{e}^{\mathrm{i}\varphi})|&<2.12^2\\
|Q_N^{(2)}(\mathrm{e}^{\mathrm{i}\varphi})|&<1.001\cdot|1-\mathrm{e}^{\mathrm{i}\varphi}|^{\tfrac{1}{2}}.
\end{align*}
Hence we obtain
\begin{align*}
\eta&<\frac{1}{4}\big(2\sqrt{3}\delta_N^{(1)}+\sqrt{3}\delta_N^{(3)}+\pi 3^{\frac{3}{4}}\delta_N^{(2)}\cdot1.001+2.12^2\pi\,3^{\frac{3}{4}}\varepsilon_N\\\nonumber
&\quad+\pi\,3^{\frac{5}{2}}\delta_N^{(2)}\varepsilon_N+3^{\frac{1}{4}}\delta_N^{(2)}\cdot1.001+2.12^2\cdot3^{\frac{1}{4}}\delta_N^{(2)}\varepsilon_N+\delta_N^{(3)}\\
&\quad+\frac{|c_1|}{4}3^{\frac{3}{4}}\delta_N^{(4)}+3\pi3^{\frac{5}{4}}\varepsilon_N+12\pi\delta_N^{(1)}+3\pi\delta_N^{(1)}+\frac{3|c_2|}{4}3^{\frac{3}{4}}\delta_N^{(4)}\big).\nonumber
\end{align*}
Thus, for $N\geq20$ we have $\eta<0.005041$. Furthermore, we can estimate
\begin{align*}
\left|\frac{h_{\mathrm{num}}^{\mathrm{approx}}\left(\mathrm{e}^{\mathrm{i}\varphi}\right)}{h_{\mathrm{denum}}\left(\mathrm{e}^{\mathrm{i}\varphi}\right)}\right|&\leq\left|\frac{h_{\mathrm{num}}^{\mathrm{approx}}\left(\mathrm{e}^{\mathrm{i}\varphi}\right)}{4|1+\mathrm{e}^{\mathrm{i}\varphi}|^{\tfrac{7}{2}}}\right|.
\end{align*} 
By construction, the right-hand side of this inequality is the absolute value of a polynomial in $1+\mathrm{e}^{\mathrm{i}\varphi}$ which takes its maximum on $[\tfrac{\pi}{3},\pi]$ at $\varphi=\pi$. For $N=20$ this maximum value is less than 0.1764. Combining both estimates we end up with $\max_{\varphi\in[\tfrac{\pi}{3},\pi]}\left|h^{(2)}(\mathrm{e}^{\mathrm{i}\varphi})\right|<0.182$.\\

\noindent Similarly, from (\ref{eqn:sqrt1-z}) we obtain at $z=1$ the series expansions
\begin{align*}
(1+z)^{\tfrac{1}{2}}&=\sqrt{2}\sqrt{1-\frac{1-z}{2}}=\sqrt{2}\left(1-\sum_{n=1}^{\infty}\frac{(2n-2)!}{n!(n-1)!2^{2n-1}}\left(\frac{1-z}{2}\right)^n\right)
\\&
=\sum_{n=0}^{\infty}\widetilde{v}_n(1-z)^n,
\\
(1+z)^{-\tfrac{1}{2}}&=\sqrt{2}\sum_{n=0}^{\infty}\frac{(2n)!}{n!n!2^{3n+1}}(1-z)^n,
\\
(1+z)^{-\tfrac{3}{2}}&=\sum_{n=0}^{\infty}\widetilde{w}_n(1-z)^n,
\end{align*}
and
\begin{align*}
\arcsin
z+\tfrac{\pi}{2}&=\pi-\sqrt{2}\sum_{n=0}^{\infty}\frac{(2n)!}{n!n!2^{3n+1}}\,\frac{1}{n+\frac{1}{2}}(1-z)^{n+\frac{1}{2}}
\\&
=\pi+\sum_{n=0}^{\infty}\widetilde{u}^{(1)}_n(1-z)^{n+\frac{1}{2}},
\end{align*}
which are convergent for $|1-z|<2$. From the latter we also obtain according series expansions
\begin{align*}
(\arcsin z+\tfrac{\pi}{2})^k&=\pi^k+\sum_{m=1}^{\infty}\widetilde{u}^{(k)}_m(1-z)^{\frac{m}{2}},\quad k=2,3,4.
\end{align*}
Let 
\begin{align*}
\widetilde{P}_N^{(k)}(z)&=\pi^k+\sum_{m=1}^{2N}\widetilde{u}^{(k)}_m(1-z)^{\frac{m}{2}},\quad
k=1,2,3,4,
\\
\widetilde{Q}_N(z)&=\sum_{n=0}^{N}\widetilde{v}_n(1-z)^{n},
\\
 \widetilde{R}_N(z)&=\sum_{n=0}^{N}\widetilde{w}_n(1-z)^{n}
\end{align*}
and
\begin{align*}
\widetilde{p}_N^{(k)}(z)&=\frac{\arcsin
  z+\tfrac{\pi}{2}-\widetilde{P}_N^{(k)}(z)}{(1-z)^{\frac{3}{2}}},\quad
k=1,2,3,4,
\\
\widetilde{q}_N(z)&=\frac{(1+z)^{\tfrac{1}{2}}-\widetilde{Q}_N(z)}{(1-z)^{\frac{3}{2}}},
\\
\widetilde{r}_N(z)&=\frac{(1+z)^{-\tfrac{3}{2}}-\widetilde{R}_N(z)}{(1-z)^{\frac{3}{2}}}.
\end{align*}
We write $h^{(2)}(z)$ as a fraction with numerator
\begin{align*}
\widetilde{h}_{\mathrm{num}}(z)&=(1+z)^{-\tfrac{3}{2}}\big(2(\arcsin
z+\tfrac{\pi}{2})(1+z)(1-z)
\\&\quad
-(\arcsin z+\tfrac{\pi}{2})^3(1+z)(1-z)\\ \nonumber
&\quad
-\pi(\arcsin z+\tfrac{\pi}{2})^2(1+z)^{\frac{3}{2}}(1-z)^{\frac{3}{2}}
\\&\quad
+z(\arcsin z+\tfrac{\pi}{2})^2(1+z)^{\frac{1}{2}}(1-z)^{\frac{1}{2}}-z^2(\arcsin z+\tfrac{\pi}{2})^3\\ \nonumber
&\quad
+\frac{c_1}{4}(\arcsin
z+\tfrac{\pi}{2})^4(1+z)^{\frac{3}{2}}-3\pi(1+z)^{\frac{5}{2}}(1-z)^{\frac{1}{2}}
\\&\quad
+4\pi(\arcsin z+\tfrac{\pi}{2})(1+z)^2(1-z)\\ \nonumber
&\quad
+\pi z(\arcsin z+\tfrac{\pi}{2})(1+z)^2-\frac{3\,c_2}{4}(\arcsin z+\tfrac{\pi}{2})^4(1+z)^{\frac{3}{2}}(1-z)\big)
\end{align*}
and denumerator
\begin{align*}
\widetilde{h}_{\mathrm{denum}}(z)&=(\arcsin z+\tfrac{\pi}{2})^4(1-z)^{\frac{3}{2}}.
\end{align*}
By $\widetilde{h}_{\mathrm{num}}^{\mathrm{approx}}(z)$ we denote the approximation of $\widetilde{h}_{\mathrm{num}}(z)$ if we replace $(\arcsin z+\tfrac{\pi}{2})^k$ by $\widetilde{P}_N^{(k)}(z)$, $\sqrt{1+z}$ by $\widetilde{Q}_N(z)$ and $(1+z)^{-\tfrac{1}{2}}$ by $\widetilde{R}_N(z)$ within $\widetilde{h}_{\mathrm{num}}(z)$. On the set $\{\mathrm{e}^{\mathrm{i}\varphi}:\varphi\in[0,\frac{\pi}{3}]\}$ we can make the following estimate for the error
\begin{align*}
\widetilde{\eta}(\varphi)&=\left|\frac{\widetilde{h}_{\mathrm{num}}\left(\mathrm{e}^{\mathrm{i}\varphi}\right)-\widetilde{h}_{\mathrm{num}}^{\mathrm{approx}}\left(\mathrm{e}^{\mathrm{i}\varphi}\right)}{\widetilde{h}_{\mathrm{denum}}\left(\mathrm{e}^{\mathrm{i}\varphi}\right)}\right|.
\end{align*}
Since on $[0,\tfrac{\pi}{3}]$ we have
\begin{align*}  
|\arcsin(\mathrm{e}^{\mathrm{i}\varphi})+\frac{\pi}{2}|&\geq|\arcsin(\mathrm{e}^{\mathrm{i}\frac{\pi}{3}})+\frac{\pi}{2}|>2.11
\end{align*}
we obtain
\begin{align*}
2.11^{4}\widetilde{\eta}(\varphi)&<
\left|h_{\mathrm{num}}(\eiphi)\right|\widetilde{r}_N(\eiphi)\\
&\quad+\bigg[2|\widetilde{p}_N^{(1)}(\eiphi))||1+\eiphi||1-\eiphi|+|\widetilde{p}_N^{(3)}(\eiphi)||1+\eiphi||1-\eiphi|\\ \nonumber
&\quad+\pi|\widetilde{P}_N^{(2)}(\eiphi)||\widetilde{q}_N(\eiphi)||1+\eiphi||1-\eiphi|^{\tfrac{3}{2}}\\
&\quad+\pi|\widetilde{p}_N^{(2)}(\eiphi)||\widetilde{Q}_N(\eiphi)||1+\eiphi||1-\eiphi|^{\tfrac{3}{2}}\\ \nonumber
&\quad+\pi|\widetilde{p}_N^{(2)}(\eiphi)||\widetilde{q}_N(\eiphi)||1+\eiphi||1-\eiphi|^{3}\\
&\quad+|\widetilde{P}_N^{(2)}(\eiphi)||\widetilde{q}_N(\eiphi)||1-\eiphi|^{\tfrac{1}{2}}\\ \nonumber
&\quad+|\widetilde{p}_N^{(2)}(\eiphi)||\widetilde{Q}_N(\eiphi)||1-\eiphi|^{\tfrac{1}{2}}\\
&\quad+|\widetilde{p}_N^{(2)}(\eiphi)||\widetilde{q}_N(\eiphi)||1-\eiphi|^{2}\\ \nonumber
&\quad+|\widetilde{p}_N^{(3)}(\eiphi)|+\frac{|c_1|}{4}|\widetilde{P}_N^{(4)}(\eiphi)||\widetilde{q}_N(\eiphi)||1+\eiphi|\\
&\quad+\frac{|c_1|}{4}|\widetilde{p}_N^{(4)}(\eiphi)||\widetilde{Q}_N(\eiphi)||1+\eiphi|\\ \nonumber
&\quad+\frac{|c_1|}{4}|\widetilde{p}_N^{(4)}(\eiphi)||\widetilde{q}_N(\eiphi)||1-\eiphi|^{\tfrac{3}{2}}|1+\eiphi|\\
&\quad+3\pi|\widetilde{q}_N(\eiphi)||1-\eiphi|^{\tfrac{1}{2}}|1+\eiphi|^2\\ \nonumber
&\quad+4\pi|1-\eiphi||1+\eiphi|^2|\widetilde{p}_N^{(1)}(\eiphi)|+\pi|1+\eiphi|^2|\widetilde{p}_N^{(1)}(\eiphi)|\\ \nonumber
&\quad+\frac{3|c_2|}{4}|\widetilde{P}_N^{(4)}(\eiphi)||\widetilde{q}_N(\eiphi)||1+\eiphi||1-\eiphi|\\
&\quad+\frac{3|c_2|}{4}|\widetilde{p}_N^{(4)}(\eiphi)||\widetilde{Q}_N(\eiphi)||1+\eiphi||1-\eiphi|\\ \nonumber
&\quad+\frac{3|c_2|}{4}|\widetilde{p}_N^{(4)}(\eiphi)||\widetilde{q}_N(\eiphi)||1+\eiphi||1-\eiphi|^{\tfrac{5}{2}}\bigg]\cdot|\widetilde{R}_N(\eiphi)|.
\end{align*}
Since $|1-\mathrm{e}^{\mathrm{i}\varphi}|$ is increasing on $[0,\pi]$ we have 
\begin{align*}
|1-\mathrm{e}^{\mathrm{i}\varphi}|^{-1}&\leq|1-\mathrm{e}^{\mathrm{i}\frac{\pi}{3}}|^{-1}=1
\end{align*}
on $[0,\frac{\pi}{3}]$. The coefficients $\widetilde{u}_n^{(1)}$, $\widetilde{v}_n$, and $\widetilde{w}_n$ are non-alternating for $n\geq1$. Furthermore, for $m\geq1$, we have $\widetilde{u}_{2m}^{(k)}>0$ and $|\widetilde{u}_{2m}^{(k)}|>|\widetilde{u}_{2m+1}^{(k)}|$, $k=2,3,4$. Since we have
\begin{align*}
|1-\mathrm{e}^{\mathrm{i}\varphi}|&\leq|1-\mathrm{e}^{\mathrm{i}\frac{\pi}{3}}|=1
\end{align*}
on $[0,\frac{\pi}{3}]$ we can conclude that there
\begin{align*}
|\widetilde{p}_N^{(k)}(\mathrm{e}^{\mathrm{i}\varphi})|&\leq|\widetilde{p}_N^{(k)}(\mathrm{e}^{\mathrm{i}\frac{\pi}{3}})|=:\widetilde{\delta}_N^{(k)},\quad k=1,2,3,4,\\
|\widetilde{q}_N(\mathrm{e}^{\mathrm{i}\varphi})|&\leq|\widetilde{q}_N(\mathrm{e}^{\mathrm{i}\frac{\pi}{3}})|=:\widetilde{\varepsilon}_N,\\
|\widetilde{r}_N(\mathrm{e}^{\mathrm{i}\varphi})|&\leq|\widetilde{r}_N(\mathrm{e}^{\mathrm{i}\frac{\pi}{3}})|=:\widetilde{\kappa}_N
\end{align*}
is satisfied. Finally, we have
\begin{align*}
|\widetilde{P}_N^{(k)}(\mathrm{e}^{\mathrm{i}\varphi})|&\leq\pi^k,\quad k=1,2,3,4,\\
|1+\mathrm{e}^{\mathrm{i}\varphi}|&\leq2,\quad |\widetilde{Q}_N(\eiphi)\leq2.\\
\end{align*}
Hence we obtain
\begin{align*}
\widetilde{\eta}&\leq\big(4\pi+2\pi^3+2\sqrt{2}\pi^3+\sqrt{2}\pi^2+\pi^3+\tfrac{|c_1|}{2}\sqrt{2}\pi^4
\\
&\qquad
+12\sqrt{2}\pi+16\pi^2+4\pi^2+\tfrac{3|c_2|}{2}\sqrt{2}\pi^4\big)\widetilde{\kappa}_N\\ \nonumber
&\quad+2\widetilde{\delta}_N^{(1)}3^{-\frac{1}{4}}+\widetilde{\delta}_N^{(3)}3^{-\frac{1}{4}}+\pi^3\widetilde{\varepsilon}_N3^{-\frac{1}{4}}+\pi\widetilde{\delta}_N^{(2)}+\pi\widetilde{\delta}_N^{(2)}\widetilde{\varepsilon}_N3^{-\frac{1}{4}}\\
&\quad+\pi^2\widetilde{\varepsilon}_N+\widetilde{\delta}_23^{-\frac{1}{2}}+\widetilde{\delta}_N^{(2)}\widetilde{\varepsilon}_N\\ \nonumber
&\quad+\widetilde{\delta}_N^{(3)}+\tfrac{|c_1|}{4}\pi^4\widetilde{\varepsilon}_N3^{-\frac{1}{4}}+\tfrac{|c_1|}{4}\widetilde{\delta}_N^{(4)}+\tfrac{|c_1|}{4}\widetilde{\delta}_N^{(4)}\widetilde{\varepsilon}_N3^{-\frac{1}{4}}+3\pi\widetilde{\varepsilon}_N\sqrt{2}+4\pi\sqrt{2}\widetilde{\delta}_N^{(1)}\\ \nonumber
&\quad+\pi\sqrt{2}\widetilde{\delta}_N^{(1)}+\tfrac{3}{4}|c_2|\pi^4\widetilde{\varepsilon}_N3^{-\frac{1}{4}}+\tfrac{3}{4}|c_2|\widetilde{\delta}_N^{(4)}+\tfrac{3}{4}|c_2|\widetilde{\delta}_N^{(4)}\widetilde{\varepsilon}_N3^{-\frac{1}{4}}.
\end{align*}
For $N\geq12$ we have $\widetilde{\eta}<0.01052$. \\

\noindent Both the functions $|\arcsin(\eiphi)+\tfrac{\pi}{2}|$ and
 $|\widetilde{h}_{\mathrm{num}}^{\mathrm{approx}}(\eiphi)(1-\eiphi)^{-{3/2}}|$
 for $N=12$ are decreasing on $[0,\tfrac{\pi}{3}]$. In order to estimate the quotient of both we divide $[0,\tfrac{\pi}{3}]$ into small subintervals of length $\tfrac{\pi}{300}$. For each subinterval it is easy to evaluate that the quotient is less than 0.1. Hence we have 
\begin{align*}
\left|\frac{\widetilde{h}_{\mathrm{num}}^{\mathrm{approx}}\left(\mathrm{e}^{\mathrm{i}\varphi}\right)}{\widetilde{h}_{\mathrm{denum}}\left(\mathrm{e}^{\mathrm{i}\varphi}\right)}\right|&<0.1
\end{align*} 
which implies together with the estimate for $\widetilde{\eta}$ that also on $[0,\tfrac{\pi}{3}]$ $\left|h^{(2)}(\mathrm{e}^{\mathrm{i}\varphi})\right|$ is less than 0.182.
\begin{figure}
\begin{center}
\includegraphics[height=5cm]{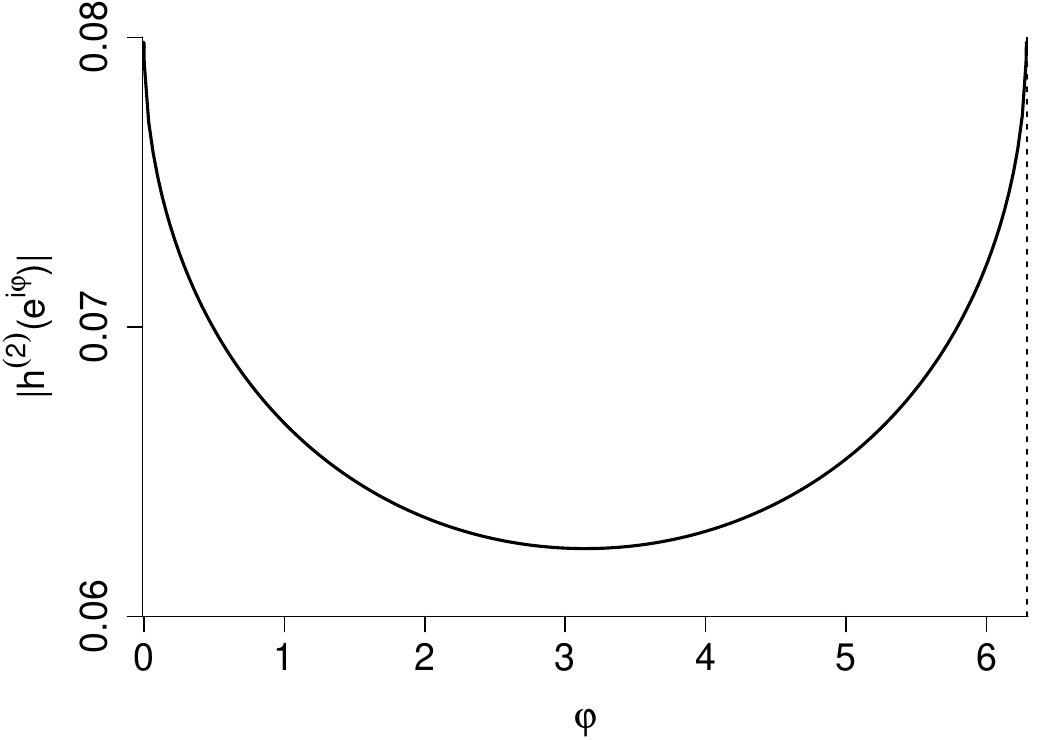}
\caption{$\left|h^{(2)}(\mathrm{e}^{\mathrm{i}\varphi})\right|$ for $\varphi\in[0,2\pi]$.}
\label{fig:Proof_abs_mon}
\end{center}
\end{figure}
\end{proof}

\subsection{Proof of Lemma \ref{lem:covExcursionSet} }
We have $P_t(1)=\mathbb{P}(Z(o)\geq t)=\Psi(t)$, and a series expansion of $P_t(\rho)$ around $\rho=1$ yields
\begin{align*}
P_t(\rho)&=\Psi(t)+\frac{\varphi(t)}{2\sqrt{\pi}}\left(-2\sqrt{1-\rho}+O((1-\rho)^{\frac{3}{2}})\right).
\end{align*} From 
this it is easy to see that in case
\begin{align*}
R'(0+)&=\lim_{\varepsilon\rightarrow0+}\frac{R(\varepsilon)-1}{\varepsilon}<0
\end{align*}
(including $R'(0+)=-\infty$) we always have
\begin{align*}
C_{\Xi_t}'(0+)&=\lim_{\varepsilon\rightarrow0+}\frac{P_t(R(\varepsilon))-P_t(1)}{R(\varepsilon)-1}\frac{R(\varepsilon)-1}{\varepsilon}\\
&=\frac{\varphi(t)}{2\sqrt{\pi}}\lim_{\varepsilon\rightarrow0+}\frac{2}{\sqrt{1-R(\varepsilon)}}\frac{R(\varepsilon)-1}{\varepsilon}=-\infty.
\end{align*}
In case $R'(0+)=0$, which means that the covariance function of $Z$ is differentiable at the origin, define
\begin{align*}
R''(0+)=2\lim_{\varepsilon\rightarrow0+}\frac{R(\varepsilon)-1}{\varepsilon^2}.
\end{align*}
Then we have
\begin{align*}
C_{\Xi_t}'(0+)&=-\frac{\varphi(t)}{\sqrt{2\pi}}\lim_{\varepsilon\rightarrow0+}\frac{\sqrt{1-R(\varepsilon)}}{\varepsilon}
=-\frac{\varphi(t)}{\sqrt{2\pi}}\sqrt{-R''(0+)}.
\end{align*}
Hence, if $R''(0+)=-\infty$, which means that the covariance function of $Z$ is not twice differentiable at the origin, then also $C_{\Xi_t}'(0+)=-\infty$. Otherwise we are in the case where $R''(0+)$ is finite and $Z$ by definition mean-square differentiable \cite{ref:Adler1981}.

\subsection{Proof of Theorem \ref{thm:covtDiff}}
By $P_t'(\rho)=\varphi(t,t,\rho)$, $E_t'(\rho)=t\varphi(t,t,\rho)+\varphi(t)\Psi(t\sqrt{(1-\rho)/(1+\rho)})$ and $C_t'(\rho)=t^2\varphi(t,t,\rho)+2t\varphi(t)\Psi(t\sqrt{(1-\rho)/(1+\rho)})+P_t(\rho)$ we obtain for $\mathrm{cov}(r)=f_t(R(r))$ according to (\ref{eqn:f_t})
\begin{equation*}
\mathrm{cov}'(r)=f_t'(R(r))R'(r),\quad f_t'(\rho)=A(\rho)P_t'(\rho)+B(\rho),
\end{equation*}
where
\begin{align*}
A(\rho)&=\frac{(t^2-\rho)P_t(\rho)^2-(2+4\rho)t\varphi(t)\Psi\bigg(t\sqrt{\frac{1-\rho}{1+\rho}}\bigg)P_t(\rho)+2(1+\rho)^2\varphi(t)^2\Psi\bigg(t\sqrt{\frac{1-\rho}{1+\rho}}\bigg)^2}{P_t(\rho)^3},\\
B(\rho)&=1+\frac{2t\varphi(t)\Psi\bigg(t\sqrt{\frac{1-\rho}{1+\rho}}\bigg)P_t(\rho)-(1-\rho^2)\varphi(t,t,\rho)^2-2(1+\rho)\varphi(t)^2\Psi\bigg(t\sqrt{\frac{1-\rho}{1+\rho}}\bigg)^2}{P_t(\rho)^2}.
\end{align*}
For each $t\in\RR$ fixed, $B(\rho)$ is bounded from below and from above for all $\rho\geq0$ since $0<\Psi(t)^2\leq P_t(\rho)\leq\Psi(t)$, $1/2\leq\Psi(t\sqrt{(1-\rho)/(1+\rho)})\leq1$ and $\varphi(t)^4\leq(1-\rho^2)\varphi(t,t,\rho)^2\leq\varphi(t)^2/(2\pi)$. Hence, if $R'(0+)=0$, we have $\mathrm{cov}'(r)=A(1)C_{\Xi_t}'(0+)$ with $C_{\Xi_t}'(0+)<0$ being finite and $A(1)=((t^2-1)\Psi(t)^2-3t\varphi(t)\Psi(t)+2\varphi(t)^2)/\Psi(t)^3$. 
With
\begin{align*}
 - 945 t^{-11}
\le \Psi(t) / \varphi(t) - (t^{-1} - t^{-3} + 3 t^{-5}-15 t^{-7} + 105 t^{-9})
\le0
, \qquad t\ge0,
\end{align*}
we have
\begin{align*}
&(t^2-1)\Psi(t)^2-3t\varphi(t)\Psi(t)+2\varphi(t)^2\\
&\;\;\;\ge \varphi^2(t)\Big[
(t^2-1) (t^{-1} - t^{-3} + 3t^{-5}-15 t^{-7} + 105  t^{-9} - 945
t^{-11} )^2 \\
&\qquad
- 3 t   (t^{-1} - t^{-3}
+ 3t^{-5}  -15 t^{-7} + 105  t^{-9}
) + 2 \Big]\\
&\;\;\;=\varphi^2(t) \bigg[
{\frac{2}{t^6}}-{\frac{30}{t^8}}-{\frac{2439}{t^{10}}}+{\frac{4935
 }{t^{12}}}-{\frac{11565}{t^{14}}}+{\frac{48195}{t^{16}}}
-{\frac{237825
 }{t^{18}}}+{\frac{1091475}{t^{20}}}-{\frac{893025}{t^{22}}}
\bigg]
\\
&\;\;\;>0,\qquad t\ge8
.
\end{align*}
Numerical inspection for $0\leq t<8$ and similar considerations for $t<0$
yield that $(t^2-1)\Psi(t)^2-3t\varphi(t)\Psi(t)+2\varphi(t)^2>0$, $t\in\RR$, that is, $A(1)$ is positive.

%
If $Z$ is not mean-square differentiable we have $R'(0+)<0$ and $C_{\Xi_t}'(0+)=-\infty$. Since $B(1)$ is bounded and $A(1)$ is positive we obtain $\mathrm{cov}'(0+)=-\infty$.
Since $R''(0+)=0$ implies that $R$
is constant \cite[Lemma~1.10.16]{ref:Sasvari1994},
Lemma \ref{lem:covExcursionSet} yields that 
$\mathrm{cov}'(0+)\in[-\infty,0)$.

For $\gamma(r)=v_t(R(r))$ according to (\ref{eqn:v_t}) we obtain
\begin{equation*}
v_t(\rho)=1-\rho+(1-\rho)^2\frac{t\varphi(t)\Psi\bigg(t\sqrt{\frac{1-\rho}{1+\rho}}\bigg)-(1+\rho)\varphi(t,t,\rho)}{P_t(\rho)},
\end{equation*}
and hence
\begin{align*}
v_t'(\rho)&=-1+\sqrt{1-\rho}\,\frac{(1+2\rho)\sqrt{1-\rho}\,\varphi(t,t,\rho)-2\sqrt{1-\rho}\,t\varphi(t)\Psi\bigg(t\sqrt{\frac{1-\rho}{1+\rho}}\bigg)}{P_t(\rho)}\\
&\quad+(1-\rho)\,\frac{(1-\rho^2)\varphi(t,t,\rho)^2-t\varphi(t)\Psi\bigg(t\sqrt{\frac{1-\rho}{1+\rho}}\bigg)(1-\rho)\varphi(t,t,\rho)}{P_t(\rho)^2}.
\end{align*}
For each $t\in\RR$ fixed, the coefficients of $\sqrt{1-\rho}$ and $1-\rho$ are bounded for all $\rho\geq0$. Thus we have $v_t'(1)=-1$ which implies $\gamma'(0+)=-R'(0+)$.

\section{Acknowledgement}
Zakhar Kabluchko has been supported by DFG-SNF FOR 916
``Statistical Regularisation and Qualitative Constraints''.



\end{document}